\documentclass{amsart}

\usepackage{color, mathtools}

\mathtoolsset{showonlyrefs}

\newcommand{\ddb}{\sqrt{-1}\partial\overline{\partial}}

\renewcommand{\[}{\begin{equation} \begin{aligned} }
\renewcommand{\]}{\end{aligned} \end{equation}}

\newtheorem{thm}{Theorem}
\newtheorem{prop}[thm]{Proposition}

\newtheorem{conj}[thm]{Conjecture}

\theoremstyle{definition}
\newtheorem{defn}[thm]{Definition}

\numberwithin{equation}{section}

\begin{document}
\title{K\"ahler-Einstein metrics}

\author{G\'{a}bor Sz\'{e}kelyhidi}
\address{Department of Mathematics\\
         University of Notre Dame\\
         Notre Dame, IN 46556}
\email{gszekely@nd.edu}

\dedicatory{Dedicated to Sir Simon Donaldson on the occasion of
  his 60th birthday}
 
\begin{abstract}
We survey the theory of K\"ahler-Einstein metrics, with particular
focus on the circle of ideas surrounding 
the Yau-Tian-Donaldson conjecture for Fano manifolds.  
\end{abstract}

\maketitle

\section{Introduction}
A starting point in the study of K\"ahler-Einstein metrics is K\"ahler's
observation~\cite{Kah33}, that for Hermitian metrics satisfying what
is now known as the
K\"ahler condition, the Einstein equations reduce to a scalar complex
Monge-Amp\`ere equation. Over the many decades since, the 
field has grown into a very rich subject 
with deep connections to nonlinear PDE, geometric
analysis, complex algebraic geometry, string theory, and others.
The goal of this survey is to give an overview of some of
these developments and in particular to showcase the
diverse ideas that have been brought to bear on the problem. 

Let us start with K\"ahler's observation, and consider a Hermitian
metric $g_{j\bar k}$ on a complex manifold $M$. The associated
(1,1)-form, or K\"ahler form, is defined to be
\[ \omega = \sqrt{-1} g_{j\bar k} dz^j\wedge
d\bar{z}^k \]
in local coordinates, and the metric $g$ is K\"ahler if
$d\omega=0$. K\"ahler showed that in this case we can locally write
the metric $g$ in terms of a potential function $\phi$:
\[ g_{j\bar k} = \frac{\partial^2\phi}{\partial z^j \partial
  \bar{z}^k}. \]
The Ricci curvature of $g$ is then given by
\[ \mathrm{Ric}_{j\bar k} = -\frac{\partial^2}{\partial z^j \partial
  \bar{z}^k} \log\det(g), \]
and so we can obtain solutions of the Einstein equation $\mathrm{Ric}
= \lambda g$, by solving the scalar equation
\[ \det\left(\frac{\partial^2 \phi}{\partial
    z^j \partial\bar{z}^k}\right) = e^{-\lambda \phi}. \]

Under certain conditions K\"ahler potentials exist globally, not just
locally.  Let us suppose that $M$ is compact.
A K\"ahler form $\omega$ on $M$ defines a cohomology class
$[\omega]\in H^2(M)$, and it is natural to consider, as
Calabi~\cite{Cal54} did, the space of all K\"ahler forms on $M$ in
a fixed cohomology class. The $\partial\bar{\partial}$-lemma states
that any other K\"ahler form $\eta\in [\omega]$ can be written as
\[ \eta = \omega + \sqrt{-1}\partial\bar{\partial} \phi \]
for a function $\phi:M\to\mathbf{R}$, and so 
the space of K\"ahler metrics in a fixed cohomology class are
parametrized by scalar functions, in analogy with a conformal class
in Riemannian geometry. 

A further important observation is that for any K\"ahler metric $g$ on
$M$, its Ricci form
\[ \mathrm{Ric}(g) = \sqrt{-1} R_{j\bar k} dz^j \wedge d\bar{z}^k \]
is a closed form in the first Chern class $c_1(M)$. 
Calabi~\cite{Cal54} conjectured that conversely any representative of
$c_1(M)$ is the Ricci form of a unique K\"ahler metric in every
K\"ahler class. This fundamental conjecture was proven by
Yau~\cite{Yau78}, by solving the complex Monge-Amp\`ere equation
\[ (\omega + \sqrt{-1}\partial\bar{\partial}\phi)^n = e^{F+c}
  \omega^n\]
for $\phi$ and a constant $c$, given a K\"ahler form $\omega$ and function $F$. 
Perhaps the most important case, which has had an enormous impact, is
when $c_1(M)=0$. In this case Yau's result implies that every K\"ahler
class on $M$ admits a unique Ricci flat metric.

More generally, if we seek a K\"ahler-Einstein metric
$\omega$ satisfying $\mathrm{Ric}(\omega) = \lambda\omega$, 
then we must have $c_1(M) = \lambda[\omega]$. In
particular if $\lambda\ne 0$, then 
either $c_1(M)$ or $-c_1(M)$ must be a K\"ahler class, and
the cohomology class $[\omega]$ is uniquely determined. When $c_1(M)$
is negative, then the works of Yau~\cite{Yau78} and Aubin~\cite{Aub78}
yield a K\"ahler-Einstein metric on $M$. It was already known
by Matsushima~\cite{Mat57}, however, that when $c_1(M)$ is positive,
i.e. $M$ is Fano, then $M$
can only admit a K\"ahler-Einstein metric if its holomorphic
automorphism group is reductive. Later Futaki~\cite{Fut83} found a different
obstruction stemming from the automorphism group, showing
that a certain numerical invariant $F(v)$ must vanish for
all holomorphic vector fields $v$ on $M$. These obstructions
rule out the existence of a
K\"ahler-Einstein metric on the blowup
$\mathrm{Bl}_p \mathbf{P}^2$ for instance. On the other hand, 
Tian~\cite{Tian90} showed that 
in the case of Fano surfaces the reductivity of the automorphism
group, or alternatively the vanishing of Futaki's obstruction, 
is actually sufficient for the existence of a K\"ahler-Einstein
metric. 

At this point let us digress briefly on parallel developments in
the theory of holomorphic vector bundles. In algebraic geometry a
basic problem is to construct moduli spaces of various objects, for
instance vector bundles over a curve. It turns out that in general it
is not possible to parametrize all vector bundles of a fixed
topological type with a nice space, but rather we need to
restrict ourselves to semistable bundles -- a notion introduced by
Mumford~\cite{MFK94}. While this is a purely
algebro-geometric notion, it was shown by
Narasimhan-Seshadri~\cite{NS65}, and later reproved by
Donaldson~\cite{Don83_1}, that stability has a differential geometric
meaning: an indecomposable vector bundle over a curve is stable if and only if
it admits a Hermitian metric with constant curvature. The
Hitchin-Kobayashi correspondence, proved by Donaldson~\cite{Don85, Don87} and
Uhlenbeck-Yau~\cite{UY86}, is the higher dimensional generalization of
this, stating that an indecomposable vector bundle is stable if and
only if it admits a Hermitian-Einstein metric. There is a particularly
rich interplay between this result for complex surfaces and
Donaldson theory~\cite{DK90} for smooth four-manifolds. 

In analogy with these results on vector bundles, Yau~\cite{Yau93}
conjectured that the existence of a K\"ahler-Einstein metric on a Fano
manifold $M$ should be related to the stability of $M$ in a suitable
sense. This conjecture was made precise by Tian~\cite{Tian97}, who
generalized Futaki's obstruction~\cite{Fut83} to the notion of
K-stability:  Tian showed
(see also Ding-Tian~\cite{DT92}) that given any
$\mathbf{C}^*$-equivariant family $\pi : X\to \mathbf{C}$
with generic fiber $X_t \cong M$, and $\mathbf{Q}$-Fano central fiber
$X_0$, the Futaki invariant $F(X)$ of the induced
vector field on the central fiber can be defined. The
$\mathbf{Q}$-Fano condition here means that $X_0$ is a normal variety
with $\mathbf{Q}$-Cartier anticanonical divisor, and this
assumption allows for a differential-geometric definition of $F(X)$.  
Tian showed moreover that if $M$ admits a
K\"ahler-Einstein metric, then $F(X)\geq 0$ for all
such families, with equality only if $X$ is a product. This 
obstruction is called K-stability, and it 
is a far reaching generalization of Futaki's obstruction. Indeed the latter can be
viewed as a special case using only product families. 

The Donaldson-Uhlenbeck-Yau theorem and Yau and Tian's conjectures on
the existence of K\"ahler-Einstein metrics can be seen as two
instances of a relationship between quotient constructions in
symplectic and algebraic geometry, due to Kempf-Ness~\cite{KN79}. This is
because in both settings the geometric structure we seek, a
Hermitian-Einstein metric or a K\"ahler-Einstein metric, can be viewed
as a zero of a moment map. This was discovered by
Atiyah-Bott~\cite{AB83} for vector bundles over curves, and
independently by Fujiki~\cite{Fuj92} and Donaldson~\cite{Don97} for
K\"ahler-Einstein metrics. In fact even more generally, constant
scalar curvature K\"ahler metrics, and the extremal
K\"ahler metrics introduced by Calabi~\cite{Cal82} fit into this
framework. 

Motivated by this, Donaldson~\cite{Don02} introduced a
generalization of K-stability for any pair $(M, L)$ of a projective
manifold $M$ equipped with an ample line bundle $L$. The definition is
similar to Tian's notion, in that we need to consider $\mathbf{C}^*$-equivariant
degenerations $\pi : X \to \mathbf{C}$ of $M$, 
compatible with the polarization $L$ of $M$. The central fiber, however, is
allowed to be a singular scheme, and the corresponding numerical
invariant, the Donaldson-Futaki invariant $DF(X)$, is defined
purely algebraically. In this generality we have

\begin{conj}[Yau-Tian-Donaldson] The manifold $M$ admits a constant
  scalar curvature K\"ahler metric in $c_1(L)$, if and only if the
  pair $(M,L)$ is K-stable. 
\end{conj}

The conjecture can be extended~\cite{GSz04} to characterize the
existence of extremal metrics, and there are
also variants for more general ``twisted'' equations by Dervan~\cite{Der14}.
One direction of the conjecture is fairly well understood, namely that
the existence of a constant scalar curvature metric implies
K-stability (see e.g. Tian~\cite{Tian97}, Donaldson~\cite{Don01,Don05} and
Stoppa~\cite{Sto08}, Berman-Darvas-Lu~\cite{BDL16}), however
the converse in general is wide open at present. 

The main subject of this survey is the case when $M$ is a Fano manifold and
$L=-K_M$, since then a constant scalar curvature metric in $c_1(L)$
is actually K\"ahler-Einstein. In this case Chen-Donaldson-Sun~\cite{CDS12, CDS13_1,
  CDS13_2, CDS13_3} proved the following breakthrough result. 

\begin{thm}\label{thm:CDS}
  A Fano manifold $M$ admits a K\"ahler-Einstein metric if and only if
  $(M,-K_M)$ is K-stable. 
\end{thm}

Our aim in this survey is not so much to describe the proof of this
result, but rather to highlight the diversity of ideas that are
in some way related to the Yau-Tian-Donaldson conjecture. There are
several other excellent surveys on the subject, such as those of
Thomas~\cite{Thomas06}, Phong-Sturm~\cite{PS08} and
Eyssidieux~\cite{Eys16}, with a focus on different aspects of the
theory. 

The solution of the conjecture in the Fano case has certainly closed a chapter, but 
it has also set the scene for a great deal of further development, 
much of which is likely yet to come. 

\subsection*{Acknowledgements} It is my pleasure to thank Simon
Donaldson for his advice and support over the years -- it would be hard
to overstate the influence that his ideas and approach to mathematics 
have had on my interests.

I also thank Julius Ross and Valentino Tosatti for helpful comments on
this survey. This work was supported in part by NSF grant DMS-1350696.

\section{The moment map picture}\label{sec:momentmap}
In this section we describe how the scalar curvature of a K\"ahler
metric can be viewed as an infinite dimensional moment map, following
Donaldson~\cite{Don97}. This point of view is invaluable in building
intuition for the problem, and in retrospect it motivates many of the
basic constructions and results that were known beforehand. We will
keep the discussion at a formal level, and not delve into the precise
definitions relating to infinite dimensional manifolds. 

Let $(X,\omega)$ denote a compact symplectic manifold, such that
$H^1(X)=0$ for simplicity. Let $\mathcal{J}$ be the space of almost
complex structures on $X$, compatible with $\omega$. The space
$\mathcal{J}$ has a natural complex structure, 
and each tangent space $T_J\mathcal{J}$ is equipped with the
$L^2$-inner product given by the metric $g_J(\cdot, \cdot) =
\omega(\cdot, J \cdot)$. This structure turns $\mathcal{J}$ into an
infinite dimensional K\"ahler manifold, and the group
$\mathcal{G} = \mathrm{Ham}(X,\omega)$ of Hamiltonian symplectomorphisms acts on
$\mathcal{J}$, preserving this K\"ahler structure. We identify the Lie
algebra of $\mathcal{G}$ with the functions $C^\infty_0(X)$
with zero mean on $X$ with respect to the volume form $\omega^n$,
through the Hamiltonian construction. We further identify
$C^\infty_0(X)$ with its dual using the $L^2$ inner product. The key calculation
is the following. 
\begin{prop}[Fujiki~\cite{Fuj92},
  Donaldson~\cite{Don97}] \label{prop:scalmomentmap}
  A moment map for the action of $\mathcal{G}$ on
  $\mathcal{J}$ is given by 
  \[ \mu : \mathcal{J} &\to C^\infty_0(X) \\
       J &\mapsto  \overline{S} - S_J, \]
  where $S_J$ is the scalar curvature of the metric $g_J$ whenever
  $J$ is integrable, and $\overline{S}$ is its average, which is
  independent of $J$. 
\end{prop}
In particular an integrable complex structure $J$ satisfies $\mu(J)=0$
if and only if the K\"ahler metric $g_J$ on $X$ has constant scalar
curvature. The precise meaning of this result is an identity relating
the linearization of the scalar curvature $S_J$ under varying the
complex structure $J$, and the infinitesimal action of Hamiltonian
symplectomorphisms on $\mathcal{J}$. Indeed, let $h\in C^\infty_0(M)$,
and let $A\in T_J\mathcal{J}$ be an infinitesimal variation of
$J$. The variation of $J$ by the Hamiltonian vector field $v_h$ is the
Lie derivative $L_{v_h}J$, and we write $DS_J(A)$ for the variation of
the scalar curvature $S_J$ in the direction $A$. The content of
Proposition~\ref{prop:scalmomentmap} is the identity
\[ \langle DS_J(A), h\rangle_{L^2} = \langle JA, L_{v_h}
  J\rangle_{L^2},  \]
which can be checked by direct calculation. See \cite{Fuj92, Don97,
  Tian00, Gau15} for the details.

Suppose for the moment that instead of the infinite dimensional group
$\mathcal{G}$ acting on $\mathcal{J}$, we had a compact group $G$
acting on a compact K\"ahler manifold $(V,\omega)$, with moment map
$\mu$. Let us assume that
$\omega$ is the curvature form of a line bundle $L\to V$ endowed with
a Hermitian metric, and so $V$ is in fact
a projective manifold. Let $G^c$ denote the complexification of $G$,
acting on $V$ by biholomorphisms. The Kempf-Ness
theorem~\cite{KN79} says, in this finite dimensional situation, that a
$G^c$-orbit contains a zero of the moment map if and only if
it is polystable. One way to define polystability, that is useful in
the infinite dimensional setting as well, is that an orbit $G^c\cdot
p$ is polystable if a suitable real valued $G$-invariant function is
proper on the orbit. More precisely we consider the function
\[ \label{eq:lognormf} f : G^c/G &\to \mathbf{R} \\
            [g] &\mapsto \log \Vert g\cdot \hat{p}\Vert, 
\]
defined using a $G$-invariant norm on $V$.
Here $\hat{p} \in L$ is a non-zero lift of $p$, and we need a lift of
the $G^c$-action to the total space of $L$ in a way compatible with
the choice of moment map $\mu$. The compatibility of the lift of the
action with the moment map $\mu$ can be expressed by the formula
\[  \label{eq:df}
df_{g\cdot p}(\sqrt{-1}\xi) = \langle \mu(g\cdot p), \xi\rangle \]
for the variation of $f$, where $g\in G^c$, and $\xi\in
\mathfrak{g}$. Note in particular that the critical points of $f$ are
precisely zeros of the moment map. 

Since the function $f$ also turns
out to be convex along geodesics in the symmetric space $G^c/G$, it is
clear that properness of $f$ corresponds to the existence of a zero of
the moment map in the orbit $G^c\cdot p$, at least if we ignore
subtleties related to the possible stabilizer of $p$. 
What is less clear, however,  is that to verify whether $f$ is proper
on $G^c / G$, it is enough to check whether $f$ is proper along each
geodesic ray in $G^c / G$ obtained from one-parameter subgroups
$\mathbf{C}^*\subset G^c$. In fact it is enough to consider only 
one-parameter subgroups of the form $t \mapsto e^{\sqrt{-1} t\xi}$ for
  $\xi\in\mathfrak{g}$ generating a circle subgroup. For such a
  one-parameter subgroup we can test the properness of $f$ by
  computing the limit
\[ \label{eq:limf'} \lim_{t\to \infty} f'(e^{\sqrt{-1}t\xi}\cdot p) = \langle
  \mu(q),\xi\rangle, \]
where $q = \lim_{t\to\infty}  e^{\sqrt{-1} t \xi}\cdot p$. Properness
of $f$ is then 
equivalent to $\langle\mu(q), \xi\rangle > 0$ whenever $q
\not\in G^c\cdot p$. This is in essence the Hilbert-Mumford numerical
criterion for stability, proved by Mumford~\cite{MFK94}, to which we refer the
reader for the detailed development of this theory. 

Let us return to the infinite dimensional setting of the action of
$\mathcal{G}$ on $\mathcal{J}$. A first issue is that the
complexification $\mathcal{G}^c$ does not exist, but we can still try
to interpret what its orbits would be if it did. Indeed in each
tangent space $T_J\mathcal{J}$ we have a subspace spanned by elements
of the form $L_{v_h}J$ giving the infinitesimal action of Hamiltonian
vector fields, and we can simply complexify this subspace. The orbits
of $\mathcal{G}^c$ then ought to be integral submanifolds of this
distribution on $\mathcal{J}$. 

Note that ultimately we are interested
in the metrics $g_J$ determined by the pairs $(\omega, J)$, and for
any diffeomorphism $f$ the metric given by $(\omega, f^*J)$ is
isometric to that given by $((f^{-1})^*\omega, J)$. We can therefore
switch our point of view from studying different complex structures on
a symplectic manifold $(X, \omega)$ to studying different K\"ahler
forms on a complex manifold $(X,J)$, as is more standard in K\"ahler
geometry. To see what this corresponds to in terms of the complexified
orbits of $\mathcal{G}$, note that when $J$ is
integrable, then $JL_{v_h}J = L_{Jv_h} J$, and at the same time we have
the formula $-L_{Jv_h}\omega =
2\sqrt{-1}\partial\overline{\partial}h$. Using the
$\partial\bar{\partial}$-lemma we can interpret this as
saying that in our infinite dimensional setting the role of the
symmetric space $G^c / G$ is played by the space of K\"ahler metrics
in the K\"ahler class $[\omega]$ (see~\cite{Don99_1}). In conclusion the Kempf-Ness
theorem suggests that the existence of a constant scalar curvature
metric in the K\"ahler class $[\omega]$ is equivalent to
stability of this class in a suitable sense. 

Let us see how this formal discussion motivates several of the basic
constructions in the field, which were actually discovered before the moment
map picture was understood:

\subsection{The Mabuchi metric}\label{sec:mabuchimetric}
In finite dimensions the metric on the symmetric space $G^c /
G$ is given by an inner product on the Lie algebra $\mathfrak{g}$. In
the infinite dimensional setting we chose the $L^2$ product on
Hamiltonian functions, which by the above discussion correspond to
variations in the K\"ahler potential. This leads to a very natural
Riemannian structure on the space of K\"ahler metrics first introduced by
Mabuchi~\cite{Mab87} and later rediscovered by Semmes~\cite{Sem92} and
Donaldson~\cite{Don99_1}. For a compact K\"ahler manifold $(M,\omega)$, 
Let us denote by
\[ \mathcal{H} = \{ \phi\in C^\infty(M)\,:\, \omega + \ddb\phi >
  0 \} \]
the space of K\"ahler potentials. For $\phi\in \mathcal{H}$ let
$\omega_\phi = \omega + \ddb\phi$ be the corresponding K\"ahler
metric. When $\omega\in c_1(L)$ for an ample line bundle $L$, then
 $\mathcal{H}$ can also be thought of as the space of positively
 curved Hermitian metrics $e^{-\phi}$ on $L$. 

Each tangent space $T_\phi \mathcal{H}$ can be
identified with $C^\infty(M)$, and the Mabuchi metric is defined by
simply taking the $L^2$ inner product: 
\[ \langle f, g\rangle_\phi = \int_M fg\, \omega_\phi^n. \]
One can show that this metric turns $\mathcal{H}$, at least
formally, into a non-positively curved symmetric space. 

Of great interest is the study of geodesics in $\mathcal{H}$. A
calculation shows that a path $\phi_t \in \mathcal{H}$ is a geodesic,
if it satisfies the equation
\[ \ddot\phi_t - \frac{1}{2}|\nabla \dot\phi_t|^2_{\omega_{\phi_t}} =
  0. \]
An important observation due to Semmes and Donaldson, however, is that
this geodesic equation is equivalent to a homogeneous complex
Monge-Amp\`ere equation. Indeed, let $A_{a,b}= S^1 \times (a,b)$ be a
cylinder, and given a path $\phi_t\in
\mathcal{H}$ for $a < t < b$, define the form
\[ \Omega = \pi^*\omega + \ddb \phi_t \]
on the product $M \times A_{a,b}$. Here $\pi:M\times A_{a,b} \to M$ is
the projection, and $\ddb$ involves the variables on $A_{a,b}$ as
well. A calculation shows that $\phi_t$ is a geodesic if and only if
$\Omega$ is non-negative and $\Omega^{n+1} = 0$, i.e. $\phi_t$ solves
the homogeneous complex Monge-Amp\`ere equation on $M\times
A_{a,b}$. When $a,b$ are finite we have geodesic segments, while if
$a$ or $b$ is infinite, then we have geodesic rays. 

Since the equation $\Omega^{n+1}=0$ is 
degenerate elliptic, the regularity theory is very
subtle. Chen~\cite{Chen00_1} showed that any two potentials $\phi_0,
\phi_1 \in \mathcal{H}$ can be connected by a unique weak
geodesic $\phi_t$, for which $\Delta \phi_t$ is bounded, using
the Laplacian on $M\times A_{a,b}$ (see also
B{\l}ocki~\cite{Bl12}). This was improved to a bound on $|\nabla^2
\phi_t|$ by Chu-Tosatti-Weinkove~\cite{CTW16} (see also
Berman~\cite{Ber14} for a weaker result
 in the projective case). It turns out that these results
are essentially optimal, since there are counterexamples to the
existence of smooth geodesics, by Lempert-Vivas~\cite{LV11} and
Darvas~\cite{Dar12} (see also Donaldson~\cite{Don02_1} and Ross-Witt
Nystr\"om~\cite{RN15}).  However even
weak geodesics are enough for many applications, as we will describe
below. 

\subsection{The K-energy}
In the finite dimensional setting we described how the existence of a
zero of the moment map is related to properness of the log-norm
functional $f$ in \eqref{eq:lognormf}. In infinite dimensions this translates to the
K-energy, defined by Mabuchi~\cite{Mab86}. The formula \eqref{eq:df}
for the variation of the log-norm functional suggests that the
K-energy $K: \mathcal{H} \to\mathbf{R}$ can be defined through its
variation. If $\phi_t\in \mathcal{H}$ is a path, then  
\[ \frac{d}{dt} K(\phi_t) &= \langle \overline{S} -
  S(\omega_{\phi_t}), \dot\phi_t\rangle_{L^2(\omega_{\phi_t})} \\
  &= \int_M \dot{\phi}_t (\overline{S} - S(\omega_{\phi_t}))\,
  \omega_{\phi_t}^n. \]
Mabuchi~\cite{Mab87} showed that the K-energy is convex along
smooth geodesics in $\mathcal{H}$, which we now see as a general
result about the log-norm functionals. It is clear from
the definition that the critical points are constant scalar curvature
metrics, and Mabuchi also showed that if two critical points $\phi_0,
\phi_1$ are connected by a smooth geodesic, then the metrics
$\omega_{\phi_0}$ and $\omega_{\phi_1}$ are isometric by an automorphism
of $M$. 

It was only much later that Berman-Berndtsson~\cite{BB14}
showed that convexity holds along weak geodesics, and as an application
proved a general uniqueness result along these lines. Note that
uniqueness in various degrees of generality has been proven previously
using other methods, see for example~\cite{BM85, Don01, Chen00_1,
  Ber13}.

In the finite dimensional setting the existence of a
critical point of the log-norm functional is equivalent to its
properness. Tian~\cite{Tian97} showed that the analogous result holds
for K\"ahler-Einstein metrics, characterizing their existence in terms
of properness of the K-energy in a suitable sense. See also
Darvas-Rubinstein~\cite{DR15} for a more precise properness statement
in the K\"ahler-Einstein case, in the presence of automorphisms. In
the general constant scalar curvature case 
Berman-Darvas-Lu~\cite{BDL16} showed one direction of this
correspondence, namely that the existence of a cscK metric implies properness of
the K-energy, as was conjectured by Tian~\cite{Tian00}. 

\subsection{The Futaki invariant}
A construction that predates both of the previous ones is Futaki's
obstruction~\cite{Fut83} to the existence of a K\"ahler-Einstein
metric on a Fano manifold $M$, analogous to the Kazdan-Warner
obstruction~\cite{KW74} for the prescribed curvature problem on the
2-sphere. 

In retrospect, Futaki's obstruction can be viewed as the first glimpse
of the obstruction to K\"ahler-Einstein metrics given by
K-stability. Recall that in the finite dimensional picture,
polystability of $p$ is related to the limit \eqref{eq:limf'} of the derivative of the
log-norm functional along the orbit $e^{\sqrt{-1}t\xi}\cdot \hat{p}$ of a
one-parameter subgroup. The simplest example is if $\xi\in \mathfrak{g}$ is in the
stabilizer of $p$, so that the one-parameter subgroup simply acts on
the line $L_p$. The quantity $\langle \mu(p), \xi\rangle$ is then the
weight of this action, and polystability requires that this weight
vanishes, since otherwise the log-norm functional would not be bounded
from below. 

The infinite dimensional analog of this weight can be defined as
follows. An element $\xi\in\mathfrak{g}$ in the stabilizer of a point
$p$ corresponds to a function $h$ on $M$, whose Hamiltonian vector
field $v_h$ preserves the complex structure of $M$ as well, i.e. $v_h$
is a holomorphic Killing field on $(M,\omega)$, and $v_h$ generates a
circle action on $M$. The corresponding weight is then
\[ \label{eq:futdef}
F(v_h) = \int_M h(\overline{S} - S(\omega))\,\omega^n. \]
Futaki used a different, but essentially equivalent definition, and showed
that $F(v_h)$ only depends on the vector field $v_h$, and not on the
metric in the class $[\omega]$ used in the formula. In addition the
invariant can be defined for any holomorphic vector field, not just
those that generate circle actions. 
If $M$ admits a constant scalar curvature metric in the class
$[\omega]$, then it is clear from the definition that $F(v)=0$ for all
holomorphic vector fields $v$. 

Tian's definition~\cite{DT92, Tian97} of K-stability is motivated
by probing the properness of the K-energy along more general families
of metrics $\phi_t\in \mathcal{H}$. From the finite dimensional
picture it is most natural to consider geodesic rays, however this is
technically rather difficult. Instead Tian used families of metrics on
$M$ obtained from embedding $M\subset \mathbf{P}^N$ into a projective
space, and then considering the restrictions of the Fubini-Study
metrics $\sigma_t^*\omega_{FS}$ pulled back under a one-parameter
family of automorphisms $\sigma_t$ of $\mathbf{P}^N$. We will discuss
this construction and others in more detail in Section~\ref{sec:Kstab}. 

\subsection{The Ding functional}
The constructions in the previous subsections apply to the general
existence problem for constant scalar curvature metrics, not just
K\"ahler-Einstein metrics. At the same time we will see that the K\"ahler-Einstein
problem has several special features. One of these is 
an alternative variational description of K\"ahler-Einstein metrics as
critical points of the Ding functional $\mathcal{D}$ defined in \cite{Ding88}. Thinking of
$\mathcal{H}$ as the space of positively curved metrics $e^{-\phi}$ on
$-K_M$, the variation of $\mathcal{D}$ along a path $\phi_t$ is
defined to be
\[ \frac{d}{dt} \mathcal{D}(\phi_t) = - \frac{1}{V}\int_M \dot\phi_t
  \omega_{\phi_t}^n + \frac{\int_M \dot\phi_t e^{-\phi_t}}{\int_M
    e^{-\phi_t}}, \]
where $V$ is the volume with respect to $\omega^n$, and we can
naturally think of $e^{-\phi_t}$ as defining volume forms on $M$. The
critical points of this functional satisfy $e^{-\phi} =
C\omega_\phi^n$, so they are K\"ahler-Einstein metrics. 

The Ding functional has many analogous properties to the K-energy,
such as the convexity along weak geodesics proved by
Berndtsson~\cite{Ber13}, but it has technical advantages over the
K-energy, since defining it requires less regularity of
$\phi$. Recently Donaldson~\cite{Don15} gave a variation of the
infinite dimensional moment map picture discussed above, in which
Berndtsson's convexity result~\cite{Ber09} gives rise to the K\"ahler
structure on $\mathcal{J}$, and the Ding functional corresponds to the
log-norm functional. The weight again recovers the Futaki
invariant, and the existence of a K\"ahler-Einstein metric is related
to properness of $\mathcal{D}$ by Tian~\cite{Tian97}.

\section{K-stability}\label{sec:Kstab}
In this section we survey the concept of K-stability of a Fano
manifold $M$, or more generally a projective manifold $M$ with an
ample line bundle $L$, from different
points of view. We first discuss the original notion for Fano
manifolds, due to Tian~\cite{Tian97}, which is fairly
differential geometric. A much more
algebro-geometric definition for general pairs $(M,L)$ 
was given by Donaldson~\cite{Don02}. In 
both of these definitions one needs to consider
$\mathbf{C}^*$-equivariant degenerations of $M$, and the main
difference is that in Tian's definition the central fiber is required
to be a $\mathbf{Q}$-Fano variety, whereas it can be an arbitrary
scheme in Donaldson's definition. This added flexibility is needed when
dealing with general polarized manifolds, but Li-Xu~\cite{LX11} showed
that in the Fano case the two notions of K-stability are
equivalent. We will now consider these two notions in more detail,
along with a more analytic approach through geodesic rays in
$\mathcal{H}$. 

\subsection{Tian's definition}
The first notion of K-stability was introduced by
Tian~\cite{Tian97}, in the context of Fano manifolds. As we have
discussed in the previous section, Tian showed that if a Fano manifold
$M$ admits a K\"ahler-Einstein metric, then the K-energy on the
K\"ahler class $c_1(M)$ is proper, and K-stability can be thought of
as probing this properness along certain families of metrics. 

Suppose that we have a $\mathbf{C}^*$-equivariant family of varieties $\pi: X \to
\mathbf{C}$, with generic fiber $\pi^{-1}(t) \cong M$ for $t\ne
0$. Assume in addition that the central fiber is normal, and
that a power of the relative anticanonical line bundle on the regular
locus extends to a relatively ample line bundle on $X$. In this
situation we call the family $X$ a special degeneration of $M$. 

The $\mathbf{C}^*$-action on such a 
 special degeneration induces a $\mathbf{C}^*$-action on 
the central fiber $X_0 = \pi^{-1}(0)$. Using that $X_0$ has relatively
mild singularities, Tian (see also Ding-Tian~\cite{DT92}) showed that
one can define the Futaki invariant of this $\mathbf{C}^*$-action on
$X_0$ using a differential geometric formula similar to
\eqref{eq:futdef}. This is then defined to be the Futaki invariant
$F(X)$ of the special degeneration $X$. Note that any
$\mathbf{C}^*$-action on the Fano manifold $M$ gives rise to a product
action on the trivial family $X= M\times \mathbf{C}$, and the Futaki
invariant of this family is simply the Futaki invariant of the
original $\mathbf{C}^*$-action. At the same time there are infinitely
many special degenerations, even if $M$ admits no
$\mathbf{C}^*$-actions. 

In order to relate special degenerations to families of metrics and
properness of the K-energy, note that any special degeneration $\pi: X
\to \mathbf{C}$ for $M$ can be
realized as a family in projective space. More precisely there is an
embedding $X \subset \mathbf{P}^N \times\mathbf{C}$, such that the
$\mathbf{C}^*$-action on $X$ is induced by the action of a
 one-parameter subgroup
$\sigma:\mathbf{C}^*\to SL(N+1)$ on $\mathbf{P}^N$. Here $M$ is
embedded in $\mathbf{P}^N\times \{1\}$ using a basis of sections of
$-rK_M$ for a suitable $r > 0$. We can now define a family of metrics
$\omega_t\in c_1(M)$ by restricting the Fubini-Study metric to
the non-zero fibers of $X$. Equivalently we can write
\[ \label{eq:testconfigpath}
\omega_t = \frac{1}{r} \sigma_{e^{-t}}^*\omega_{FS}|_M. \]
Ding-Tian~\cite{DT92} showed that with suitable normalizing factors
which we omit,
\[ \label{eq:limddtK} \lim_{t\to\infty} \frac{d}{dt} K(\omega_t) = F(X), \]
as suggested by \eqref{eq:limf'} in the finite dimensional picture.
Note, however, that the family $\omega_t$ is usually not a geodesic ray. 

With these results in mind we have the following definition, due to
Tian~\cite{Tian97}. 
\begin{defn} \label{defn:TianKstab}
  A Fano manifold $M$ is K-stable, if $F(X)\geq 0$ for all special
  degenerations $X$ of $M$, with equality only for product
  degenerations. 
\end{defn}
In the same paper, Tian showed that if $M$ admits a K\"ahler-Einstein
metric then the K-energy is proper, and as a consequence he showed the
following fundamental result (see also Berman~\cite{Ber12} for the
case when $M$ admits holomorphic vector fields). 
\begin{thm} If a Fano manifold admits a K\"ahler-Einstein metric, then
  it is K-stable. 
\end{thm}

\subsection{The Donaldson-Futaki invariant}
The central fiber of a special degeneration has fairly mild
singularities, and so a differential geometric definition of the
Futaki invariant was possible. On the other hand
Donaldson~\cite{Don02} assigned an invariant, called the
Donaldson-Futaki invariant, to essentially arbitrary
$\mathbf{C}^*$-equivariant degenerations of a polarized manifold
$(M,L)$, by giving a purely algebro-geometric definition of the
Futaki invariant for a $\mathbf{C}^*$-action on any polarized scheme
$(V,L)$. 

To give the definition, note that a $\mathbf{C}^*$-action $\lambda$ on $(V,L)$
induces actions on the spaces of sections $H^0(V,kL)$, and so in
particular for each $k$ we have a total weight $w_k$. For large $k$ we
have expansions
\[ \label{eq:Hilbpoly}
\dim H^0(V,kL) &= a_0 k^n + a_1 k^{n-1} + \ldots \\
    w_k &= b_0 k^{n+1} + b_1 k^n + \ldots, \]
and the Donaldson-Futaki invariant of $\lambda$ is defined to be
\[DF(\lambda) = \frac{a_1}{a_0}b_0 -  b_1.\] 
When $V$ is smooth, the equivariant Riemann-Roch formula can be used
to show that this coincides with Futaki's differential geometric
definition. 

Such polarized schemes with $\mathbf{C}^*$-actions arise naturally as
the central fibers of test-configurations. 

\begin{defn}
  Let $(M,L)$ be a polarized manifold. A test-configuration for
  $(M,L)$ with exponent $r$ is a $\mathbf{C}^*$-equivariant flat
  family $\pi : (X,\mathcal{L})\to \mathbf{C}$, such that
  $\mathcal{L}$ is relatively ample, and
  \[ (\pi^{-1}(t), \mathcal{L}|_{\pi^{-1}(t)}) \cong (M, rL), \]
  for any $t\ne 0$. In addition it is natural to require that 
  the total space $X$ is normal (see Li-Xu~\cite{LX11} and
  Ross-Thomas~\cite{RT04}). 

  The Donaldson-Futaki invariant $DF(X,\mathcal{L})$ of the
  test-configuration is defined to be the Donaldson-Futaki invariant
  of the induced $\mathbf{C}^*$-action on the central fiber.  
\end{defn}

Given this definition, K-stability can be defined as follows, in
analogy with Definition~\ref{defn:TianKstab}.

\begin{defn}\label{defn:Kstab}
A polarized manifold $(M,L)$ is K-stable, if
$DF(X,\mathcal{L})\geq 0$ for all test-configurations for
$(M,L)$, with equality only if $X\cong M\times \mathbf{C}$. 
\end{defn}

Using this definition, the following result was shown by
Stoppa~\cite{Sto08}, building on work by Donaldson~\cite{Don05} and
Arezzo-Pacard~\cite{AP06}.
\begin{thm}
  Suppose that $M$ admits a constant scalar curvature metric in
  $c_1(L)$, and it has no nonzero holomorphic vector fields. Then
  $(M,L)$ is K-stable. 
\end{thm}

The result can be extended to the case when $M$ has holomorphic vector
fields, and also to extremal metrics (see \cite{SSz09},
\cite{BDL16}). As we stated in the introduction, the
Yau-Tian-Donaldson conjecture is the converse of this result, saying
that if $(M,L)$ is K-stable, then there is a constant scalar curvature
metric in $c_1(L)$. However it is likely that
actually a stronger notion of stability is needed in general, in view of examples
of Apostolov-Calderbank-Gauduchon-T\o{}nnesen-Friedman~\cite{ACGT3},
that are shown to be unstable in a suitable sense by
Dervan~\cite{Der16}. One possibility for such a stronger stability
notion is provided by the formalism of filtrations~\cite{GSz13, Ny10},
while another is the concept of uniform K-stability~\cite{Der14,BHJ15}. 

In order to compare Definitions~\ref{defn:TianKstab} and \ref{defn:Kstab},
let us point out that it is fairly easy to construct interesting
test-configurations, which are not special degenerations, using for instance
deformation to the normal cone. This was explored in detail by
Ross-Thomas~\cite{RT04, RT06}. One can simply take any subscheme
$Z\subset M$, and let $X = \mathrm{Bl}_{Z\times \{0\}} M\times
\mathbf{C}$, with a suitable relatively ample line bundle
$\mathcal{L}$. For instance when $Z$ is a smooth submanifold of $M$,
then the central fiber of $X$ will be isomorphic to $\mathrm{Bl}_Z M\cup P_Z$, where
$P_Z = \mathbf{P}(N_{Z}\oplus \underline{\mathbf{C}})$ is the
projective completion of the normal bundle of $Z$ in $M$, and $P_Z$ 
is glued along its zero section to the blowup $\mathrm{Bl}_Z M$ along
its exceptional divisor. In fact Odaka~\cite{Oda13_1} showed that by
blowing up ``flag ideals'' of $M\times\mathbf{C}$ instead of just
subschemes, one can essentially recover all test-configurations, and
using this approach Odaka-Sano~\cite{OS12} and Dervan~\cite{Der15}
were able to prove the K-stability of certain varieties. 

With this in mind, it appears that in the Fano case
Definition~\ref{defn:Kstab}  is more
restrictive than Tian's Definition~\ref{defn:TianKstab}, 
since test-configurations are much more
general than special degenerations. It is quite remarkable then that
for Fano manifolds the two notions turn out to be equivalent. This was
first proven by Li-Xu~\cite{LX11} purely algebro-geometrically, using the
minimal model program. Roughly speaking the minimal model program
allowed them to modify an arbitrary test-configuration into a special
degeneration, while controlling the sign of the Donaldson-Futaki
invariant at each step.  A more differential geometric proof also
follows from Chen-Donaldson-Sun's proof~\cite{CDS12} of the
YTD-conjecture for Fano manifolds.

One suggestive example is to consider a polarized toric manifold $(M,L)$,
with Delzant polytope $P$. It is natural in this case to only allow
torus equivariant test-configurations. The
only torus equivariant test-configurations with normal central fiber
are product configurations
induced by a $\mathbf{C}^*$-action on $M$ and indeed, when $M$ is Fano,
then Wang-Zhu~\cite{WZ04} showed that the only obstruction to the
existence of a K\"ahler-Einstein metric is that given by the Futaki
invariants of these $\mathbf{C}^*$-actions. On the other hand, as
shown by Donaldson~\cite{Don02}, any rational piecewise linear convex
function on $P$ gives rise to a test-configuration for $(M,L)$ and
there are (non-Fano) examples where these give an obstruction to the
existence of a cscK metric, not detected by product configurations.

\subsection{Intersection theoretic formula}
An alternative formula for the Donaldson-Futaki invariant in terms of
intersection products has been very useful in more algebro-geometric
developments. It was observed by Wang~\cite{Wang12}, and it is 
also related to the CM-polarization of
Tian~\cite{Tian97, PT06_1}. To explain it, note that any test-configuration
$(X,\mathcal{L})$ can be extended trivially at infinity to obtain a
$\mathbf{C}^*$-equivariant family 
\[ (\overline{X}, \overline{\mathcal{L}}) \to \mathbf{P}^1. \]
The line bundle $\overline{\mathcal{L}}$ is relatively ample, and by
taking the tensor product with a line bundle pulled back from
$\mathbf{P}^1$ we can assume that it is actually ample. A calculation
shows that in terms of this family 
the Donaldson-Futaki invariant of a test-configuration of exponent $r$ is
\[ \label{eq:DFintersect}
DF(X, \mathcal{L}) = \frac{n}{n+1}  \mu(M,rL)\,
  (\overline{\mathcal{L}})^{n+1} +
  \overline{\mathcal{L}}^n.K_{\overline{X}/\mathbf{P}^1},\]
using the intersection product on $\overline{X}$. 
Here for a polarized variety  $(M,L)$, the ``slope'' $\mu(M,L)$ is defined by
\[ \mu(M,L) = \frac{-K_M.L^{n-1}}{L^n}, \]
and up to a constant multiple is it the average scalar curvature
$\overline{S}$ of a K\"ahler metric in $c_1(L)$. 

This reformulation of the Donaldson-Futaki invariant has various
advantages, as shown for instance in the works of Li-Xu~\cite{LX11},
Odaka~\cite{Oda13} and others. Here we just mention one,
namely the extension of K-stability to non-algebraic K\"ahler
manifolds due to Dervan-Ross~\cite{DR16} and Sj\"ostr\"om
Dyrefelt~\cite{SD16}. While the expansions \eqref{eq:Hilbpoly} do not
make sense in the absence of a line bundle, one can make sense of the
intersection product \eqref{eq:DFintersect} even in the K\"ahler
case.

\subsection{Geodesic rays}\label{sec:geodesicrays}
In Section~\ref{sec:momentmap} we described how in the finite
dimensional moment map picture, stability can be tested using geodesic
rays in the symmetric space $G^c / G$. Donaldson~\cite{Don99_1}
formulated conjectures saying that in an analogous way
 geodesic rays in $\mathcal{H}$ can
be used to detect the existence of a constant
scalar curvature metric. Since that time there has been enormous
progress on our understanding of geodesic rays, although
these conjectures are still mostly open except in the Fano case.

We have seen in
\eqref{eq:testconfigpath} that a special degeneration or a
test-configuration for $M$ gives rise to a path in the space
$\mathcal{H}$ of K\"ahler potentials. Unless we have a product
test-configuration, this path cannot be expected to be a geodesic ray
in $\mathcal{H}$, but rather it is a geodesic in a finite
dimensional space of \emph{Bergman metrics}, i.e. those obtained by
restricting the Fubini-Study metric. It turns out that the relation
\eqref{eq:limddtK} between the Futaki invariant of a
test-configuration and the asymptotic derivative of the K-energy
along the corresponding Bergman geodesic does not hold for general
test-configurations. The general formula for the limit has been
obtained by Paul~\cite{Paul12} in terms of hyperdiscriminant and Chow
polytopes, leading to an alternative notion of stability. 

To relate this to geodesic rays in $\mathcal{H}$, note that
a given test-configuration $X$ for $M$ can be realized as a family in
projective spaces of arbitrarily large dimension, and in this way we
obtain not one, but a whole sequence of Bergman geodesics of metrics
on $M$ from $X$ using the formula
\eqref{eq:testconfigpath}. Phong-Sturm~\cite{PS06} showed that one can
pass to a limit, and obtain a geodesic ray in $\mathcal{H}$ in a
suitable weak sense, with an arbitrary initial point $\phi_0$. 

One can
also directly construct such a weak geodesic ray in $\mathcal{H}$ from
the test-configuration $X$ in the following way (see \cite{AT03}, \cite{CT08}, \cite{PS10},
\cite{Ber12} for this in various degrees of generality). Let us denote
by $X_\Delta$ the family $X$ restricted to the unit disk
$\Delta\subset\mathbf{C}$. We use the ``initial point'' $\phi_0$ to define a
metric $e^{-\phi_0}$ on the line bundle $\mathcal{L}$ over $\partial
X_\Delta$. The geodesic ray is then obtained by finding
 an $S^1$-invariant metric $e^{-\phi}$ on $\mathcal{L}$ over
$X_\Delta$, with positive curvature current, and solving the
homogeneous complex Monge-Amp\`ere equation $(\ddb \phi)^{n+1} = 0$ on
the interior of $X_\Delta$, in the sense of pluripotential
theory. The existence of such a solution, and its regularity
properties are discussed by Phong-Sturm~\cite{PS10}. Over the punctured
disk $\Delta^*$ the family $X$ is biholomorphic to $M\times \Delta^*$,
and so the metric $e^{-\phi}$ on $\mathcal{L}$ induces a family of
metrics on $L\to M$, which is the geodesic ray we were after.

There has been a lot of work relating the behavior of the K-energy
along such a geodesic ray to the Donaldson-Futaki invariant of the
test-configuration (see e.g. \cite{CT08}, \cite{PS10}). 
The first sharp result in this direction is due to
Berman~\cite{Ber12}, using the Ding functional instead of the
K-energy, in the case when $M$ is Fano. 
He shows that along a geodesic ray $\phi_t$ constructed from
a test-configuration $(X,\mathcal{L})$ as above, we have
\[ \label{eq:DFformula2}  DF(X,\mathcal{L})= \lim_{t\to\infty} \frac{d}{dt}
  \mathcal{D}(\phi_t) + q, \]
where $q\geq 0$ is a rational number determined by the central fiber
of the test-configuration. An analogous formula for the K-energy, for
test-configurations with smooth total space, was subsequently obtained
by Sj\"ostr\"om Dyrefelt~\cite{SD16}. The relation~\eqref{eq:DFformula2}
 led Fujita~\cite{Fuj15_1} to
study the notion of Ding stability, where the Donaldson-Futaki
invariant is replaced by the asymptotic derivative of the Ding
functional. As an important application, he showed
that projective space has the maximal volume amongst K\"ahler-Einstein
manifolds of a fixed 
dimension, using that K\"ahler-Einstein manifolds are stable in this
sense. This direction has been taken much further in recent work of
Fujita~\cite{Fu16}, Li~\cite{Li15}, Li-Xu~\cite{LX16, LX17_1},
Liu-Xu~\cite{LX17}, leading to new examples of K-stable varieties,
as we mention in Section~\ref{sec:examples}. 

A different direction was pursued by
Boucksom-Hisamoto-Jonsson~\cite{BHJ15}, who introduced the point of
view of thinking of a test-configuration for $(M,L)$ as a non-Archimedean metric
on $L$. From this point of view the asymptotic derivative of the Ding
functional can be seen as a non-Archimedean Ding functional. 
Building on these ideas, and using techniques from non-Archimedean
analysis~\cite{BFJ08}, as well as advances on the geometry of the
space of metrics $\mathcal{H}$ due to Darvas~\cite{Dar15},
Berman-Boucksom-Jonsson~\cite{BBJ15} proposed a variational approach to
proving Theorem~\ref{thm:CDS}. This approach relies far less on
differential geometric methods than Chen-Donaldson-Sun's
approach~\cite{CDS12}, and as such it may eventually lead to
a solution of the Yau-Tian-Donaldson conjecture for singular Fano
varieties. At present, however, even in the smooth case it yields a
weaker result than~\cite{CDS12} since it requires assuming uniform
K-stability. 

\subsection{Finite dimensional approximation}\label{sec:finitedim}
In the previous section we mentioned the idea of realizing a geodesic
in $\mathcal{H}$ as a limit of finite dimensional Bergman
geodesics. This idea of approximating the geometry of $\mathcal{H}$
with that of larger and larger spaces of Bergman metrics is a
fundamental one, going back at least to work of
Tian~\cite{Tian90_1} on the problem of
approximating an arbitrary K\"ahler metric $\omega\in c_1(L)$ on $M$ with
a sequence of metrics 
\[ \omega_k = \frac{1}{k} \phi_k^*\omega_{FS}. \]
Here $\phi_k : M \to \mathbf{P}^{N_k}$ are projective
 embeddings given by bases of sections of $kL$, for sufficiently large
 $k$. Fixing a Hermitian
metric on $L$ with curvature form $\omega$, we have induced $L^2$
inner products on each space of sections $H^0(kL)$. It is then natural
to use an orthonormal basis $\{s^{(k)}_0,\ldots, s^{(k)}_{N_k}\}$ 
of $H^0(kL)$ to define the embedding $\phi_k$. A calculation shows
that
\[ \label{eq:phikomega}
 \frac{1}{k} \phi_k^*\omega_{FS} = \omega + \ddb \log
\rho_{\omega,k}, \]
where the function
\[ \rho_{\omega,k} = \sum_{i=0}^{N_k} |s^{(k)}_i|^2 \]
is called the Bergman kernel. Tian gave an asymptotic expansion of
this function as $k\to\infty$, which was later refined by
Zelditch~\cite{Zel98}, Lu~\cite{Lu98}, Catlin~\cite{Cat99} and others,
showing that
\[ \label{eq:TYZexp}
\rho_k = 1 + \frac{S(\omega)}{2} k^{-1} + O(k^{-2}), \text{ as
}k\to\infty. \]
The lower order terms can also be computed, and are all given in terms
of covariant derivatives of the curvature tensor of $\omega$ (see
Lu~\cite{Lu98}). Note that normalising constants in the coefficients depend on
conventions for defining the $L^2$ product on $H^0(kL)$. In addition,
integrating the expansion over $M$ recovers the
Hirzebruch-Riemann-Roch formula. 
Using this expansion in \eqref{eq:phikomega} we see that $\omega_k -
\omega = O(k^{-2})$, and so indeed we can approximate an arbitrary
metric $\omega$ with Bergman metrics. 

Donaldson~\cite{Don01, Don04} took this point of view further, realizing
that it is advantageous to approximate the entire moment map package
described in Section~\ref{sec:momentmap} by a sequence of analogous
finite dimensional problems, a process that can be thought of as
quantization~\cite{Don01_1}. A starting point is the work of
Zhang~\cite{Zh96} and 
Luo~\cite{Luo98}, relating the GIT stability of the Chow point of a
projective submanifold $V\subset \mathbf{P}^N$ under the action of
$SL(N+1)$ to the existence of a special embedding of $V$. Letting $n$
denote the dimension of $V$, we define
the matrix
\[ M(V)_{jk} = \sqrt{-1}\int_V \frac{Z_j \overline{Z}_k}{\sum_{i=0}^{N}
  |Z_i|^2}\, \omega_{FS}^n, \]
where the $Z_i$ are the homogeneous coordinates on $\mathbf{P}^N$. Luo
shows that GIT stability is related to the existence of \emph{balanced
  embeddings}, for which the matrix $M(V)$ is a multiple of the
identity matrix. In fact $M(V)$, or rather its projection to
$\mathfrak{su}(N+1)$ can be viewed as a moment map, and zeroes of this
moment map are the balanced embeddings. 

Given any embedding $V\subset \mathbf{P}^N$, we can consider the
restriction of the Fubini-Study metric on $V$, and then construct a
new embedding using an $L^2$ orthonormal basis of sections of
$O(1)|_V$ as above. It turns out that $V\subset
\mathbf{P}^N$ is a balanced embedding if this new embedding coincides
with the old one (up to the action of $SU(N+1)$), i.e. if the Bergman
kernel is constant. The expansion \eqref{eq:TYZexp} then suggests that
for a given pair $(M,L)$, balanced embeddings using a basis of
$H^0(kL)$ for very large $k$ ought to be related to constant scalar
curvature metrics in $c_1(L)$. Indeed, Donaldson~\cite{Don01}
 shows that if a cscK
metric exists, and $M$ admits no holomorphic vector fields, then $M$
also admits balanced embeddings using a basis of $H^0(kL)$ for sufficiently large
$k$. Since such balanced embeddings are unique up to the action of
$SU(N+1)$, this implies in particular that if a cscK metric exists, it
is unique. Moreover
$M$ satisfies an asymptotic GIT stability condition, which implies
K-semistability.

By ``quantizing''
the K-energy as well, Donaldson~\cite{Don04} showed that if
a cscK metric exists, and $M$ has no holomorphic vector fields, then
the K-energy is bounded below on $\mathcal{H}$, generalizing
Bando-Mabuchi's result~\cite{BM85} in the K\"ahler-Einstein
case. Going further, using such finite dimensional approximations, 
Donaldson~\cite{Don05} showed that a
destabilizing test-configuration leads to a nontrivial lower bound for
the Calabi functional: for any metric $\omega\in c_1(L)$, and any
test-configuration $(X,\mathcal{L})$ for $(M,L)$ we have 
\[ \Vert S(\omega) - \overline{S}\Vert_{L^2} \geq \frac{
  -DF(X,\mathcal{L})}{\Vert (X,\mathcal{L})\Vert_2}, \]
where a suitable $L^2$-type norm on test-configurations is used for
normalizing the Donaldson-Futaki invariant. This result
shows, in particular, that the existence of
a cscK metric implies K-semistability, even when $M$ has holomorphic
vector fields. 

These ideas have since been developed further in many papers. A very
small sample of such developments is 
\cite{Mab04, Don09_1, SZ10, Fine10, RZ10, Li11_2, Fine12,  CS12, ST15, Sey16}.

\section{Geometric and algebraic limits}
Many of the ideas that we have discussed so far apply not only to the
problem of K\"ahler-Einstein metrics, but also to more general
constant scalar curvature metrics. In this section this will no
longer be the case, since it will be crucial that the Einstein
equation controls the Ricci curvature of the metric. Our goal is to
present some of the ideas involved in Chen-Donaldson-Sun's
proof of Theorem~\ref{thm:CDS} \cite{CDS12}, that is 
the existence part of the
Yau-Tian-Donaldson Conjecture in the Fano case.

Suppose that we have a Fano manifold $M$, and we are trying to find a
metric $\omega\in c_1(M)$ satisfying $\mathrm{Ric}(\omega) =
\omega$. There are several natural ``continuity methods'' for
deforming a given initial metric $\omega_0\in c_1(M)$ to a K\"ahler-Einstein
metric. One approach is to choose a reference metric
$\alpha\in c_1(M)$, and
then try to find a family of metrics $\omega_t$ for $t\in [0,1]$ satisfying
\[ \label{eq:contmethod}
 \mathrm{Ric}(\omega_t) = t\omega_t + (1-t)\alpha. \]
Using Yau's theorem~\cite{Yau78} we can solve this equation for $t=0$, and
using the implicit function theorem Aubin~\cite{Aub84} showed that the
set of $t$ for which a solution exists is open. The remaining
difficulty is then to understand what can happen to a sequence of
solutions $\omega_{t_k}$ as $t_k \to T$, for some $T\in (0,1]$. 
This approach was used in
\cite{DSz15}, based on the techniques of Chen-Donaldson-Sun. 

Chen-Donaldson-Sun~\cite{CDS12} used a variant of this,
proposed by Donaldson~\cite{Don09}, where the $\omega_t$ are
K\"ahler-Einstein metrics on the complement $M\setminus D$ of a
suitable smooth divisor $D\subset M$, while along $D$ the $\omega_t$
have singularities modeled on a two-dimensional cone with cone angle
$2\pi t$. In analogy with \eqref{eq:contmethod}, the metrics
$\omega_t$ satisfy
\[ \label{eq:conicalcontinuity} \mathrm{Ric}(\omega_t) = t\omega_t + (1-t)[D], \]
where $[D]$ is the current of integration along $D$. Here the Ricci
curvature is defined as the curvature current of a metric on $-K_M$
induced by the volume form of $\omega_t$. The advantage of this
continuity path over \eqref{eq:contmethod} is that on $M\setminus D$
the metric $\omega_t$ is Einstein, at the expense of introducing
singularities along $D$. At the same time in algebraic geometry it is
natural to consider pairs $(M,D)$, and metrics with cone
singularities along $D$ are useful differential geometric counterparts
(see e.g. \cite{CGP13}). The relevant linear theory for the implicit
function theorem was developed by Donaldson~\cite{Don12}, with further
refinements by Jeffres-Mazzeo-Rubinstein~\cite{JMR16}, and as in the
case of \eqref{eq:contmethod} one needs to understand the limiting
behavior of a sequence $\omega_{t_k}$ as $t_k\to T$. 

A somewhat different strategy is given by the Ricci flow, originally
introduced by Hamilton~\cite{Ham82}. On a Fano manifold $M$ the
(normalized) K\"ahler-Ricci flow is given by the parabolic equation
\[ \frac{\partial}{\partial t} \omega_t = \omega_t -
  \mathrm{Ric}(\omega_t), \]
with initial metric $\omega_0\in c_1(M)$. It was known for a long time
that the flow exists for all time~\cite{Cao85}, and the main
difficulty was understanding the behavior of $\omega_t$ as
$t\to\infty$. Building on Perelman's deep results on the Ricci
flow~\cite{Per02} as well as on ideas in Chen-Donaldson-Sun's proof,
Chen-Wang~\cite{CW14} developed the necessary convergence theory for
understanding this limiting behavior (see also Bamler~\cite{Bam16}). 
This was then used by
Chen-Sun-Wang~\cite{CSW15} to give a proof of Theorem~\ref{thm:CDS}
using the Ricci flow. 

We now return to the continuity method \eqref{eq:contmethod}, and a
sequence $\omega_{t_k}$ with $t_k\to T$. 
If such a sequence (or a
subsequence) converges to a solution $\omega_T$ of \eqref{eq:contmethod}, then either
we can further increase $t$ using the implicit function theorem, or
$T=1$ and we obtain the K\"ahler-Einstein metric that we seek. The key
question is then to understand what happens when the sequence
$\omega_{t_k}$ diverges. In this case we need to construct a special
degeneration for $M$ with non-positive Futaki invariant. 

In broad outline the strategy is the following. We first define a
certain limit space $Z$ out
of the sequence $(M, \omega_{t_k})$ (up to choosing a subsequence),
and show that this limit can be given the structure of an
algebraic variety, and in addition $Z$ can also be viewed as an
algebro-geometric limit of the images of a sequence of projective embeddings
$\phi_k : M\to \mathbf{P}^N$. For simplicity suppose that $T = 1$. Then the limit
$Z$ is shown to admit a possibly singular K\"ahler-Einstein metric. If
$Z$ is biholomorphic to $M$, then we obtain the K\"ahler-Einstein
metric that we were after, but if not, then there is a special
degeneration for $M$ with central fiber $Z$. This necessarily has
vanishing Futaki invariant, and so $M$ is not K-stable. The argument
is a little more involved when $T < 1$, since then the limit $Z$
admits a ``twisted'' K\"ahler-Einstein metric, and we obtain a special
degeneration with a strictly negative Futaki invariant. 

In the next two sections we give some more details on the construction
of the limit $Z$, and on relating the geometric and algebraic limits. It
will be helpful to consider a more general sequence $(M_k,\omega_k)$
than that obtained from the equations~\eqref{eq:contmethod}, even
allowing the manifolds $M_k$ to vary. Some aspects of the problem,
such as the convergence theory, are easier when we consider a sequence
of K\"ahler-Einstein manifolds, however others, such as using the
K-stability assumption, are easier along a continuity method such as
\eqref{eq:contmethod} (see for instance Donaldson~\cite{Don11} for the difficulty
in using a sequence of K\"ahler-Einstein metrics).

\subsection{Gromov-Hausdorff limits}\label{sec:GHlim}
In this section we will consider a sequence $(M_k, \omega_k)$ that are
either all K\"ahler-Einstein with volumes controlled from below, or
they arise from \eqref{eq:contmethod} or another similar continuity
method. 
From a geometric point of view the advantage of working with
such metrics as opposed to those with just constant
scalar curvature is that control of the Ricci curvature allows us to
extract a Gromov-Hausdorff limit of the sequence $(M_k,
\omega_k)$. This is a consequence of the Gromov compactness
theorem~\cite{Gro07} and the Bishop-Gromov volume comparison
theorem~\cite{CE08}. For simplicity let us denote by $M_k$ the metric
space $(M_k, \omega_k)$. Up to choosing a subsequence of the $M_k$, the
limit is a metric space $(Z,d)$, such that we have metrics 
(distance functions) $d_k$ on the disjoint unions $M_k\sqcup Z$,
extending those on $M_k$ and $Z$,
 satisfying the following: for all $\epsilon > 0$, the $\epsilon$-neighborhoods of
both $M_k$ and $Z$ cover all of $M_k\sqcup Z$ for sufficiently large
$k$. 

While at first $(Z,d)$ is just a metric space, a series of works by
Cheeger-Colding~\cite{CC97, CC00, CC00_1} developed a detailed
structure theory (see also Cheeger~\cite{Ch01} for an exposition). The
general theory applies to Riemannian manifolds 
with only a lower bound on the Ricci curvature,
and also allows collapsing phenomena, but our sequence $M_k$ satisfies the
non-collapsing condition 
\[ \mathrm{Vol}(B(q_k,1))  >c > 0, \]
for a constant $c > 0$, where $q_k \in M_k$ are some basepoints. This
again uses the Bishop-Gromov inequality, together with Myers's Theorem
and the fact that we control the volumes of the $M_k$. 

The
basic concept is that of a tangent cone. For any $p\in Z$,
and a sequence $r_k \to \infty$ the Gromov compactness theorem can be
applied to the sequence of pointed manifolds $(Z, r_kd, p)$ and, up to
choosing a subsequence, we obtain a limit metric space $Z_p$, called a
tangent cone of $Z$ at $p$. A fundamental
result of Cheeger-Colding is that in our situation these tangent cones are
metric cones. More precisely for each tangent cone $Z_p$ there is a
length space $Y$ of diameter at most $\pi$, such that $Z_p$ is the
completion of $Y \times (0,\infty)$ using the metric
\[ d\big( (y_1, r_1), (y_2,r_2)\big) = r_1^2 + r_2^2 - 2r_1r_2 \cos
d_Y(y_1,y_2). \]
We write this space as $C(Y)$. 

A point $p\in Z$ is called regular, if a
tangent cone $Z_p$ is isometric to $\mathbf{R}^{m}$ (in our situation
$m=2n$). In fact in this case each tangent cone
at $p$ is Euclidean. This means that if we take a sufficiently small
ball centered at $p$, and scale it up to unit size, then it will be
very close to the Euclidean ball in the Gromov-Hausdorff sense. But
recall that for sufficiently large $k$, the metric space $M_k$ is
very close to $Z$. This implies that we can find points $p_k\in M_k$,
and small balls centered at $p_k$ which, scaled to unit size, are
close to Euclidean in the Gromov-Hausdorff sense. If  the metrics
$\omega_k$ were Einstein, then results of Anderson~\cite{An90} and
Colding~\cite{Col97}  would show that
we can choose holomorphic coordinates
centered at $p_k$ with respect to which the components of $\omega_k$
are controlled in $C^{2,\alpha}$. If all points of $Z$
were regular, this would mean that there is a uniform scale at which
we have $C^{2,\alpha}$ control of the metrics $\omega_k$. In
particular we could take a limit of the K\"ahler structures 
$(M, \omega_k)$, and turn $Z$ into a smooth K\"ahler
manifold. Note, however, that even in this case $Z$ may not be
biholomorphic to $M$.

When the $\omega_k$ are not actually Einstein, but rather arise from
\eqref{eq:contmethod}, then
Anderson's result still applies at the regular points~\cite{Sz13_1}. 
However,  in general not all points of $Z$ are
regular, and instead Cheeger-Colding
defined a stratification based on how close each tangent cone is to
being Euclidean. More precisely, for $0\leq k \leq m-1$, the set $S_k$
is defined to be those points $p\in Z$ for which no tangent cone $Z_p$
splits off an isometric factor of $\mathbf{R}^{k+1}$. They then showed
that the Hausdorff dimension of $S_k$ is at most $k$ and in addition 
$S_{m-1} = S_{m-2}$, so $Z$ is regular outside a
codimension 2 set. A key technical step in Chen-Donaldson-Sun's work 
 is to improve this dimension estimate of the singular set from the Hausdorff dimension to
the Minkowski dimension - i.e. to show that there is a constant $K$
for which the singular set $S_{m-1}$ can be covered by $Kr^{2-m}$
balls of radius $r$, for any $r\in (0,1)$. For even better estimates on
the Minkowski dimension of the singular set for Einstein manifolds, 
see Cheeger-Naber~\cite{CN11, CN14}.

This discussion  applies also to the tangent cones $Z_p$, as well as to
further iterated tangent cones. In our situation this implies that
each tangent cone  $Z_p$ is a Ricci flat K\"ahler cone, at least outside of a
codimension 2 set. In addition on the regular set in $Z_p$, the Ricci
flat K\"ahler metric is given by $\ddb r^2$, where $r$ is the distance
from the vertex of $Z_p$. 

\subsection{The partial $C^0$-estimate}
The theory described in the previous section gives a good
understanding of the ``geometric'' limit of the spaces $(M_k,
\omega_k)$. In this section we incorporate the holomorphic structure
through the partial $C^0$-estimate, which is an idea due to
Tian~\cite{Tian90} used in his work on K\"ahler-Einstein metrics on
Fano surfaces. For simplicity we will assume here that $M_k=M$, and
the sequence $\omega_k$ is obtained from \eqref{eq:contmethod}.

Suppose that $-lK_M$ is very ample for some large $l > 0$. The metrics
$\omega_k$ define
Hermitian metrics $h_k$ on $-lK_M$, with curvature forms
$l\omega_k$. As in Section~\ref{sec:finitedim}, we have 
a natural projective embedding $\phi_k : M\to \mathbf{P}^{N}$
obtained by using an orthonormal basis $\{s_i\}$ of $H^0(M, -lK_M)$ for the
$L^2$ inner product
\[ \langle s,t\rangle_{L^2} = \int_M \langle s, t\rangle_{h_k}\,
(l\omega_k)^n, \]
defined by $h_k$. As before, the metric $\omega_k$ is related to the
 pullback of the Fubini-Study metric under
$\phi_k$ by the Bergman kernel $\rho_{\omega_k, l}$. 

In Section~\ref{sec:finitedim} we saw the expansion~\eqref{eq:TYZexp}
of $\rho_{\omega_k,l}$ as
$l\to\infty$ for a fixed metric $\omega_k$, whereas here the main
interest is in obtaining bounds for $\rho_{\omega_k, l}$ for
sufficiently large $l$, which apply uniformly to a sequence of metrics
$\omega_k$. 
Tian~\cite{Tian90} conjectured that if, as in our situation, we have a
positive lower bound for the Ricci curvatures of $\omega_k$, and a
lower bound for the volumes, then for sufficiently large $l$ we have a
bound $\rho_{\omega_k, l} > d > 0$ from below. This is called the
partial $C^0$-estimate, and it is equivalent to
showing that for any point $p\in M$, the bundle 
$-lK_M$ has a holomorphic section $s$ satisfying $\Vert
s\Vert_{L^2}=1$, and $|s|_{h_k}^2(p) > d$. Note that if $-lK_M$ is very
ample, then $\rho_{\omega, l}$ is positive for any metric $\omega\in
c_1(M)$, and the key point is to have a uniform lower bound along the
sequence of metrics $\omega_k$. 

For the details on how to use the partial $C^0$-estimate to define the
structure of an algebraic variety on $Z$, see Tian~\cite{Tian90} and
more generally Donaldson-Sun~\cite{DS12}. Here let us just note that a
key consequence of the estimate is that there is a uniform bound on
the derivatives of the embeddings $\phi_k : (M, \omega_k)
\to\mathbf{P}^N$, independent of $k$. Up to further increasing the
multiple of $-K_M$ that we use, this eventually leads to the result
that the Gromov-Hausdorff limit $Z$ is homeomorphic to the
algebro-geometric limit of the images $\phi_k(M)$ in projective
space. 

Let us say a few words on the proof of the
partial $C^0$-estimate. Fixing a point $p\in M$, for each metric
$\omega_k$ we must construct
holomorphic sections $s_k$ of $-lK_M$, for sufficiently large $l$, such
that $\Vert s_k\Vert_{L^2, h_k} = 1$, and $|s_k|^2_{h_k}(p) > d > 0$
for some $d$ independent of $k$. The simplest situation is when the
geometry of each $\omega_k$ is very well controlled near $p$ at a
suitable scale. For
instance this would be the case if $\omega_k$ did not actually depend
on $k$ (or more generally if the $\omega_k$ were to converge to a smooth
metric $\omega$ on $M$). In this case, at a sufficiently small scale,
the metrics $\omega_k$ appear to be very close to the Euclidean
space $(\mathbf{C}^n, \omega_{Euc})$.  The trivial line bundle over
$\mathbf{C}^n$ has metric $e^{-l|z|^2}$ with curvature
$l\omega_{Euc}$, and so for sufficiently large $l$, there will be a
holomorphic section $s$ with $|s|^2(0)=1$ and $\Vert s\Vert_{L^2}=1$
but with $s$ decaying rapidly away from the origin. Using cutoff
functions this section can be glued onto $M$ to produce a smooth
section $s_k$ of $-lK_M$ satisfying 
\[ |s_k|^2_{h_k}(p)\sim 1,\quad  \Vert s_k\Vert^2_{L^2,h_k}\sim 1,
  \text{ and }\Vert \overline{\partial} s_k\Vert_{L^2, h_k} \sim 0. \]
Using the H\"ormander $L^2$-technique, these sections can be perturbed
slightly to obtain the required holomorphic sections of $-lK_M$. 
This method was used by Tian~\cite{Tian90_1} to find a weak form of the asymptotic 
expansion \eqref{eq:TYZexp} of the Bergman kernel as $l\to\infty$, for a fixed metric
$\omega$. As we have described in Section~\ref{sec:finitedim}, this
expansion has been very  influential in the study of constant scalar
curvature K\"ahler metrics. 

In general we do not have such good uniform control of the local geometry of
the metrics $\omega_k$. The first more general situation that was
understood was the case of Fano surfaces, by Tian~\cite{Tian90}. He
proved an orbifold compactness theorem for K\"ahler-Einstein surfaces,
using techniques employed in the study of Yang-Mills
connections by Uhlenbeck~\cite{Uhl82_1, Uhl82_2} (see also
Anderson~\cite{And89} and Bando-Kasue-Nakajima~\cite{BKN89}). In
effect this implies that at a suitable scale, the local geometry of a
non-collapsed sequence of K\"ahler-Einstein surfaces
$(M_k, \omega_k)$ is modeled either on flat $\mathbf{C}^2$,
or its quotient by a finite group. Using this together with the $L^2$
technique, Tian showed the partial $C^0$-estimate for such a sequence
of K\"ahler-Einstein surfaces. 

There was little further progress until the work of
Donaldson-Sun~\cite{DS12}, who proved
the partial $C^0$-estimate for  non-collapsing sequences
 of K\"ahler-Einstein manifolds $(M_k, \omega_k)$, generalizing the
 result of Tian. In the proof they combined the $L^2$ technique with the
Cheeger-Colding structure theory that we described above. Very
roughly the structure theory, as described above, implies that the
local geometry of the $(M_k, \omega_k)$ near a point $p_k$
is modeled by Ricci flat cones
$C(Y)$. Although $Y$ may have singularities, the estimate on the size
of the singular set ensures that we can find a cutoff function $\eta$,
supported on the regular part of $C(Y)$, for which $\Vert \nabla
\eta\Vert_{L^2}\ll 1$. Taking $l$ sufficiently large, the line
bundle $-lK_{M_k}$ with its Hermitian metric $h_k$ is modeled, at
least on the regular part, by the
trivial bundle over $C(Y)$, with trivializing section $s$ whose norm
decays exponentially fast. One can then use $\eta s$ to obtain a smooth
section of $-lK_{M_k}$ which is approximately holomorphic in an
$L^2$-sense, and which can then be perturbed to a genuine holomorphic
section $s_k$. A gradient estimate is then used to show that
$|s_k|^2(p_k) > d$ for a controlled constant $d > 0$, while $\Vert
s_k\Vert_{L^2}\sim 1$, as required. 

Shortly after Donaldson-Sun's work, the partial $C^0$-estimate was
proven in various other settings. In their proof of the
Yau-Tian-Donaldson conjecture, Chen-Donaldson-Sun~\cite{CDS13_2,
  CDS13_3} proved it for K\"ahler-Einstein metrics with conical
singularities, 
while Phong-Song-Sturm~\cite{PSS12} extended Donaldson-Sun's work to
K\"ahler-Ricci solitons using also some techniques of
Tian-Zhang~\cite{TZ12}. In~\cite{Sz13_1} we showed that it holds along
the continuity method \eqref{eq:contmethod}. 
 The ideas have also been applied to the
K\"ahler-Ricci flow, by Tian-Zhang~\cite{TZ13} in dimensions up to 3,
and Chen-Wang~\cite{CW14} in general. These results on the K\"ahler-Ricci flow
also led to a proof of Tian's general conjecture on the partial
$C^0$-estimate, for a non-collapsed sequence of K\"ahler manifolds
with only a positive lower bound on the Ricci curvature (see
Jiang~\cite{Jiang13} for dimensions up to 3, and Chen-Wang~\cite{CW14}
in general).

\subsection{The K-stability condition}\label{sec:singular}
In the previous sections we have described how out of a sequence
$(M,\omega_k)$ along the continuity method \eqref{eq:contmethod}, or
from a sequence of K\"ahler-Einstein metrics, or even from a sequence
along a solution of the K\"ahler-Ricci flow one can construct a
limit $Z$, which has the structure of an algebraic variety. This limit
variety can be singular, although it is always normal with log
terminal singularities. In addition one can 
show that $Z$ admits a singular K\"ahler metric, which
satisfies either an analog of the equation \eqref{eq:contmethod} with
$t=T$ (where recall that $T= \lim_{k\to\infty} t_k$),  is
K\"ahler-Einstein, or is a K\"ahler-Ricci soliton in the case of the
Ricci flow. In this section we will explain how the K-stability
assumption on
$M$ can be used to show that in the limit we obtain a
K\"ahler-Einstein metric on $M$. In addition we discuss some of the
developments regarding singular K\"ahler-Einstein metrics, as
they play an important role in the proof of the Yau-Tian-Donaldson conjecture. 

Let us suppose that we are working with a
continuity method such as \eqref{eq:contmethod}, we have
$T=1$, and for simplicity assume that $Z$ is smooth. Note that $Z$ then
admits a K\"ahler-Einstein metric, and the goal is to show that if $M$
is K-stable, then actually $Z\cong M$.  To see this, note first that
Matsushima's Theorem~\cite{Mat57}  implies that the
automorphism group of $Z$ is reductive. On the other hand, $Z$ can
also be viewed as an algebraic limit of the projective varieties
$\phi_k(M) \subset \mathbf{P}^N$, for a sequence of embeddings
$\phi_k: M\to\mathbf{P}^N$, and so in a suitable Hilbert scheme, $Z$
can be viewed as an element in the closure of the $GL(N+1)$-orbit of
$M$. From this point of view it is the stabilizer of $Z$ under the
$GL(N+1)$ action that is reductive, and
then an appropriate version of the Luna slice theorem~\cite{Lun73,
  Don10} implies that there is a special degeneration for $M$, whose
central fiber is $Z$. Futaki's result~\cite{Fut83} furthermore implies
that the Futaki invariant of this special degeneration must vanish. If
$M$ is K-stable, then necessarily $Z\cong M$, and so $M$
admits a K\"ahler-Einstein metric. 

Suppose now that $T < 1$, and that for simplicity $Z$ is
smooth. Under Chen-Donaldson-Sun's continuity method
\eqref{eq:conicalcontinuity} one obtains a limiting K\"ahler-Einstein
metric on $Z$, with cone singularities along a divisor $D_\infty$.
The strategy is then to
apply the argument used above, but to the pairs $(M, D)$ and $(Z,
D_\infty)$. More precisely, one shows, using the equation on $Z$, 
that in a suitable Hilbert scheme of pairs the stabilizer of
$(Z,D_\infty)$ is reductive, and from this one constructs a
test-configuration $X$ for $(M,D)$ with central fiber $(Z,D_\infty)$. 
The Donaldson-Futaki
invariant has a natural generalization to pairs~\cite{Don12}, called
the log Futaki invariant $DF(X, (1-T)D)$, depending on the divisor $D$
as well as the ``angle'' parameter $T$ such that $DF(X, 0D)=DF(X)$. It
is shown in \cite{CDS13_3} that $DF(X, (1-T)D)=0$, since in fact the log
Futaki invariant of any vector field on the pair $(Z, (1-T)D_\infty)$
vanishes. This uses the existence of a conical K\"ahler-Einstein
metric on $Z$. 
In addition $DF(X,D) > 0$, which corresponds to an existence result for
K\"ahler-Einstein metrics with cone singularities that have small cone
angles. Since the log Futaki invariant is linear in the angle
parameter, it follows that $DF(X) < 0$. 

Using the continuity method \eqref{eq:contmethod} is a bit more
cumbersome, because in that case instead of working with pairs of the
form $(M,D)$ that can be parametrized by a Hilbert scheme, we need to
work with pairs $(M, \alpha)$, where $\alpha$ is a (1,1)-form, or even
a current, and the Luna slice theorem does not apply directly in this
infinite dimensional setting. In \cite{DSz15} we overcome this problem
by an approximation argument, which we now outline. 

The starting point
is to show that out of a sequence of solutions to
\eqref{eq:contmethod} with $t_i \to T$, 
 we obtain a limiting metric $\omega_T$ on the space $Z$ satisfying
 the equation 
\[ \mathrm{Ric}(\omega_T) = T\omega_T + (1-T) \beta, \]
where $\beta$ is a positive $(1,1)$-current on $Z$. This can be used
to show that the automorphism group of the pair $(Z, \beta)$, suitably
defined, is reductive, and moreover a version of the log
Futaki invariant of the pair $(Z, (1-T)\beta)$ vanishes for all vector
fields. The missing piece for exploiting this as above 
is to find a test-configuration for $M$ with
central fiber $Z$, under which the form $\alpha$ on $M$ converges to
the current $\beta$ on $Z$, but in the end we are not able to do this. 
Instead, in order to apply a result such
as Luna's slice theorem, the idea is to approximate the currents
$\alpha$ and $\beta$ with currents of integration along divisors.
 
Note first that according to
Shiffman-Zelditch~\cite{SZ99}, we can write $\alpha$ as an average of
currents of integration $[M\cap H]$ for hyperplane sections of $M$,
under a suitable projective embedding $M\subset \mathbf{P}^N$. As
discussed above, the geometric limit space $Z$ is also an algebraic
limit $\lim_{k\to\infty} \rho_k(M)$ for automorphisms $\rho_k$ of
$\mathbf{P}^N$. We can arrange that the limits $\rho_\infty(H) = \lim_{k\to \infty}
\rho_k(H)$ exist for all hyperplanes $H$ and so $\beta$ is an average
of the hyperplane sections $[Z \cap \rho_\infty(H)]$ of $Z$. This can
be used to show that a suitable tuple $(Z, \rho_\infty(H_1), \ldots,
\rho_\infty(H_l))$ has the same automorphism group as $(Z,\beta)$. One
simple example to keep in mind is when $Z$ is a line in
$\mathbf{P}^2$, and $\beta$ is the Fubini-Study metric. In this case
the automorphism group of $(Z,\beta)$ is trivial under our definition,
and so we must ``mark'' $Z$ with at least 3 hyperplane
sections to have a tuple with the same automorphism group. 

We can argue as before with these tuples, and we obtain a
test-configuration $X$ for $M$ with central fiber $Z$, such that in
addition the limits of the hyperplanes $H_1,\ldots, H_l$ under this
test-configuration are $\rho_\infty(H_1),\ldots,
\rho_\infty(H_l)$. However, since we do not know that the form
$\alpha$ converges to $\beta$ under this test-configuration, we cannot
say anything about the Donaldson-Futaki invariant of $X$. 
To overcome this we choose enough hyperplane
sections so that by an approximation argument the log Futaki invariant
of the tuple is close to that of the pair $(Z, \beta)$. It remains
then to use this to show that $DF(X) < 0$, contradicting K-stability
of $M$, unless we actually have $Z\cong M$ and $\beta = \alpha$, so
that the continuity method can be continued past $t=T$. 

To conclude this section we make some remarks about the situation when 
the limit space $Z$ is not smooth. The overall strategy remains the
same, but there are substantial technical difficulties in carrying it
through. To start with, one needs to make sense of the
K\"ahler-Einstein equation (or more general equations such as
\eqref{eq:contmethod}) on a singular variety $Z$. Working on a
resolution of singularities $\pi: X\to Z$, the equation is equivalent
to a complex Monge-Amp\`ere equation
\[ (\omega + \ddb \phi)^n = e^{F-\phi} dV, \]
on $X$, where $dV$ is a smooth volume form, however the function $F$
need not be bounded, and the form $\omega$ is only non-negative. More
precisely we have $0 \leq e^F \in L^p$ for some $ p > 1$, and $\omega$ is
also big, i.e. $\int_X \omega^n > 0$. It cannot be expected that
$\phi$ is twice differentiable, and so the equation needs to be
defined in a weak sense. This was done for the local theory of the
complex Monge-Amp\`ere operator by Bedford-Taylor~\cite{BT76}. Further
important progress in understanding such Monge-Amp\`ere equations with
right hand side in $L^p$ was made by Ko\l{}odziej~\cite{Kol98} using
techniques of pluripotential theory. A detailed study of
K\"ahler-Einstein metrics on singular varieties is given in
Eyssidieux-Guedj-Zeriahi~\cite{EGZ09}.

The lack of regularity of such a K\"ahler-Einstein metric at the
singular points means that results such a Matsushima's Theorem on the
automorphism group and the vanishing of the Futaki invariant cannot be
generalized in a straightforward way. Instead
Chen-Donaldson-Sun~\cite{CDS13_3} proved these results in the singular
setting using the theory of geodesics and Berndtsson's convexity and uniqueness
results~\cite{Ber13}. 

It is natural to expect that Theorem~\ref{thm:CDS} can be extended to
$\mathbf{Q}$-Fano varieties, and indeed Berman~\cite{Ber12} showed
that $\mathbf{Q}$-Fano varieties admitting (singular)
K\"ahler-Einstein metrics are K-stable. The converse, however, is so
far only known for smoothable $\mathbf{Q}$-Fanos, by work of
Spotti-Sun-Yao~\cite{SSY14}.

\section{Applications}
In this final section we briefly survey some applications of the
Yau-Tian-Donaldson conjecture for Fano manifolds, and the techniques
involved in its solution. The first is perhaps the most direct type of
application, namely obtaining new examples of K\"ahler-Einstein
manifolds. The second and third are more theoretical, on understanding the
moduli space of Fano K\"ahler-Einstein manifolds, and the
behavior of singular K\"ahler-Einstein metrics at the singular
points.

\subsection{New examples}\label{sec:examples}
According to Theorem~\ref{thm:CDS}, K-stability is a necessary and
sufficient condition for a Fano manifold $M$ to admit a
K\"ahler-Einstein metric. As such it is natural to expect that the
result can be used to show the existence of new K\"ahler-Einstein
metrics. Unfortunately for any given Fano manifold $M$ of dimension at
least two, there are infinitely many possible special degenerations,
and so at present there is no general method for testing
K-stability. On the other hand there are special circumstances in
which one can show that certain varieties admit a K\"ahler-Einstein
metric. 

In the two-dimensional case, Tian~\cite{Tian90} used the
$\alpha$-invariant~\cite{Tian87} and its refinements, to show that any
Fano surface with reductive automorphism group admits a
K\"ahler-Einstein metric. More generally, calculations of the
$\alpha$-invariant, and related concepts such as Nadel's multiplier
ideal sheaves~\cite{Nad90}, lead to many examples of K\"ahler-Einstein
manifolds (see e.g. \cite{CP02, DK01}). Recently Fujita~\cite{Fu16} showed that
Fano $n$-folds $M$ with $\alpha$-invariant $\alpha(M)=n/(n+1)$, which
are borderline for Tian's criterion~\cite{Tian87}, are still K-stable,
and as such they admit K\"ahler-Einstein metrics. 

A different source of examples where K-stability can be checked
effectively is Fano manifolds $M$ with large automorphism groups, since in
\cite{DSz15} we showed that it is enough to check special
degenerations that are compatible with the automorphism group. For
instance if $M$ is a toric manifold, then the only torus equivariant
special degenerations for $M$ are products, and so $M$ admits a
K\"ahler-Einstein metric whenever the Futaki invariants of all vector
fields on $M$ vanish (this result is originally due to
Wang-Zhu~\cite{WZ04}, using different methods). This result has
recently been generalized to reductive group compactifications by
Delcroix~\cite{Del15} using analytic methods. Subsequently, by
classifying equivariant special degenerations, Delcroix~\cite{Del16}
generalized this result to spherical varieties, giving new examples of
K\"ahler-Einstein manifolds. 

A more general setting is complexity-one spherical varieties~\cite{Tim11}, for
instance Fano $n$-folds with an effective action of an
$(n-1)$-dimensional torus. The equivariant special degenerations of these
complexity-one T-varieties were
classified by Ilten-S\"u\ss~\cite{IS15}, using the theory of
polyhedral divisors~\cite{AH06} and as an application they obtained
several new examples of K-stable Fano 3-folds. 

Theorem~\ref{thm:CDS} was generalized to the setting of K\"ahler cones
in \cite{CSz2} and, combined with the methods of
Ilten-S\"u\ss~\cite{IS15}, this lead to new examples of Ricci
flat K\"ahler cones. As we discussed in Section~\ref{sec:GHlim}, such
K\"ahler cones arise as tangent cones of Gromov-Hausdorff limits of
K\"ahler-Einstein manifolds, but they have also been studied
extensively from the point of view of Sasakian geometry~\cite{GMSY,
  BG08, Sp09}. To give a simple example, consider the affine variety
$M_{p,q} \subset \mathbf{C}^4$, given by the equation $xy + z^p +
w^q=0$. This has a natural action of a 2-torus, acting diagonally on
$\mathbf{C}^4$.  Using K-stability, it is shown in \cite{CSz2} that
$M_{p,q}$ admits a Ricci flat K\"ahler cone metric compatible with the
torus action, if and only if $2p > q$ and $2q > p$. These two
conditions correspond to the calculation of a 
Donaldson-Futaki invariant for two torus equivariant degenerations of
$M_{p,q}$, to the hypersurfaces given by the equations $xy + z^p=0$
and $xy + w^q=0$ respectively.

For manifolds without such large symmetry groups, the general approach of
Tian~\cite{Tian90} for Fano surfaces has recently been revisited 
in the work of Odaka-Spotti-Sun~\cite{OSS16},
Spotti-Sun~\cite{SS17} and Liu-Xu~\cite{LX17}. Here, given a family of
Fano manifolds some elements of which are known to be
K\"ahler-Einstein, one tries to find new K\"ahler-Einstein manifolds
using a continuity method in the family. The key problem is to
classify the possible Gromov-Hausdorff limits of K\"ahler-Einstein
manifolds in the family. Using this method, and combined with work of
Fujita~\cite{Fu16}, Liu-Xu~\cite{LX17} showed 
that all smooth 3-dimensional Fano hypersurfaces admit
K\"ahler-Einstein metrics. Ultimately such a study can lead to a
complete understanding of the K\"ahler-Einstein moduli space and its
compactification, as we will discuss in the next section. 

While a general algorithm for testing K-stability seems to be some
ways off at present, it seems likely that these techniques
will lead to many more examples. 

\subsection{The moduli space of K\"ahler-Einstein manifolds}
A basic question in algebraic geometry is to understand the moduli
space of a certain manifold, as well as its compactifications. Perhaps
the most well known example is the Deligne-Mumford compactification of
the moduli space of hyperbolic Riemann surfaces. This compactification
of the moduli space of smooth Riemann surfaces is obtained by adding
certain ``stable'' curves with nodal singularities, and the degeneration
of smooth curves to these singular ones can in turn be modeled
differential geometrically by degenerating families of hyperbolic
metrics. 

Canonically polarized manifolds are a higher dimensional
generalization of hyperbolic Riemann surfaces and compact moduli
spaces are obtained by adding the ``KSBA-stable'' varieties (see 
Koll\'ar-Shepherd-Barron~\cite{KSB88} Alexeev~\cite{Ale96}, and the
survey~\cite{Kol13}). It is interesting to note that the KSBA-stable
varieties were subsequently found to coincide with the K-stable
varieties by work of Odaka~\cite{Oda13}, and furthermore
Berman-Guenancia~\cite{BG14} showed that these all
admit (singular) K\"ahler-Einstein metrics. To complete the analogy
with hyperbolic Riemann surfaces it remains to understand the
convergence of smooth K\"ahler-Einstein metrics to the singular ones
on a metric level. We emphasize that in the canonically polarized
setting all smooth manifolds admit K\"ahler-Einstein metrics with
negative Ricci curvature by Yau's theorem~\cite{Yau78}, and
K-stability only enters in picking out the right type of singular
varieties to allow. 

By contrast, in the Fano case one cannot form a Hausdorff moduli space
containing all the smooth manifolds, since there are non-trivial families over
$\mathbf{C}$ that are isomorphic to a product over $\mathbf{C}^*$. For
instance any non-trivial test-configuration for a Fano manifold with smooth
central fiber gives such a family. One can hope, however, that
K-stable Fano manifolds form good moduli spaces, and that a
compactification is obtained by including certain singular K-stable
varieties. Moreover this compact moduli space should agree with the
space of K\"ahler-Einstein Fano manifolds, compactified using the
Gromov-Hausdorff topology. A  first step is to understand the
neighborhood of a K\"ahler-Einstein manifold in the moduli
space. In~\cite{GSz08} (see also Br\"onnle~\cite{BrThesis}) we proved
that given a K\"ahler-Einstein manifold $M$, its small deformations
that admit KE metrics are precisely the ones that are (poly)stable for
the action of $\mathrm{Aut}(M)$ on the deformation
space of $M$. Thus a small neighborhood of $M$ in the moduli space is
modeled by a GIT quotient for the action of $\mathrm{Aut}(M)$.  

In the case of surfaces, a global study of the moduli space of
K\"ahler-Einstein manifolds was
achieved by Odaka-Spotti-Sun~\cite{OSS16}, who explicitly identified the
compactified moduli space in terms of certain algebro-geometric 
compactifications. In higher dimensions Spotti-Sun-Yao~\cite{SSY14}
and Li-Wang-Xu~\cite{LWX16}
showed that smoothable $K$-stable $\mathbf{Q}$-Fano varieties,
i.e. the ones that we expect to appear in the compactified moduli
space, all admit singular K\"ahler-Einstein metrics, which in addition
can be obtained as Gromov-Hausdorff limits of smooth KE
metrics. Regarding the global structure of the moduli space, the
smoothable K-stable $\mathbf{Q}$-Fano varieties are parametrized by a
proper algebraic space by work of Odaka~\cite{Oda15} and
Li-Wang-Xu~\cite{LWX16}, and moreover the K-stable smooth Fanos form a
quasi-projective variety by Li-Wang-Xu~\cite{LWX15}. It remains to
understand more precisely what varieties these moduli spaces
parametrize in analogy with Odaka-Spotti-Sun's work~\cite{OSS16} on
surfaces, but as we mentioned above, there are special cases where
progress has been made, such as in the case of degree 4 del Pezzo
manifolds by Spotti-Sun~\cite{SS17} and cubic threefolds by
Liu-Xu~\cite{LX17}.

\subsection{Asymptotics of singular K\"ahler-Einstein metrics}
In Section~\ref{sec:singular} we discussed how singular
K\"ahler-Einstein metrics play an important role in the proof of
Theorem~\ref{thm:CDS}. Of particular interest are those singular
K\"ahler-Einstein varieties, which arise as limits of smooth
K\"ahler-Einstein manifolds, in view of the importance of forming
complete moduli spaces of Fano manifolds~\cite{OSS16, SSY14, Oda15,
  LWX16}, as well as Calabi-Yau manifolds~\cite{Zhang16}. 
We do not yet have a detailed
understanding of the behavior of such singular K\"ahler-Einstein
metrics near the singular points, although much progress has been
made. Such an understanding would be necessary in
order to use differential geometric techniques on these singular
spaces (see e.g. \cite{CGP13, GT16} for such applications). 

A first step in understanding the metric behavior near the singular
points is to understand the metric tangent cone, as discussed in
Section~\ref{sec:GHlim}. In the general Riemannian context the tangent
cone may depend on the sequence of rescalings chosen to construct it, but in
the setting of K\"ahler-Einstein manifolds Donaldson-Sun~\cite{DS15}
showed that the tangent cone is unique,  and it has the structure of
an affine variety. They also made progress towards an
algebro-geometric description of the metric tangent cone, and
conjectured that in general the tangent cone at a point
$x\in M$ of a singular K\"ahler-Einstein metric on $M$ can be
determined from the germ of the singularity $(M,x)$ using
K-stability. Note that it is not always the case that the germ of the tangent
cone at the vertex is biholomorphic to the germ $(M,x)$. For example
(see~\cite{Sz17}, \cite{HNpreprint}) there is a Calabi-Yau metric on a
neighborhood of the origin in the hypersurface $x_0^k + x_1^2 + x_2^2
+ x_3^2=0$ in $\mathbf{C}^4$ whose tangent cone at the origin is given
by $\mathbf{C}\times \mathbf{C}^2 / \mathbf{Z}_2$, if $k > 4$. 

An important special case of this conjecture for singular Calabi-Yau
varieties was proven by Hein-Sun~\cite{HS16}, namely when the
germ $(M,x)$ is isomorphic to the germ $(C,0)$ of a Ricci flat K\"ahler cone
$C$ at its vertex, and $C$ satisfies some additional technical
conditions. This uses some ideas of Li~\cite{Li15} and
Li-Liu~\cite{LL16} on an alternative characterization of
K-semistability based on the work of Fujita~\cite{Fuj15_1} we mentioned
in Section~\ref{sec:geodesicrays}. 
Using this, Hein-Sun showed that the singular Ricci flat
metric on $M$ near $x$ is asymptotic, in a strong sense, to the Ricci flat metric on
$C$. Recently, using more algebro-geometric techniques,
Li-Xu~\cite{LX17_1} have also made significant progress towards
resolving Donaldson-Sun's conjecture, and they have announced a
complete solution.

\bibliographystyle{acm} \bibliography{../mybib}

\begin{thebibliography}{100}

\bibitem{Ale96}
{\sc Alexeev, V.}
\newblock Log canonical singularities and complete moduli of stable pairs.
\newblock {\em arXiv:9608013\/}.

\bibitem{AH06}
{\sc Altmann, K., and Hausen, J.}
\newblock {Polyhedral divisors and algebraic torus actions}.
\newblock {\em Math. Ann. 334\/} (2006), 557--607.

\bibitem{And89}
{\sc Anderson, M.}
\newblock {Ricci curvature bounds and Einstein metrics on compact manifolds}.
\newblock {\em J. Amer. Math. Soc. 2}, 3 (1989), 455--490.

\bibitem{An90}
{\sc Anderson, M.}
\newblock Convergence and rigidity of manifolds under {R}icci curvature bounds.
\newblock {\em Invent. Math. 97\/} (1990), 429--445.

\bibitem{ACGT3}
{\sc Apostolov, V., Calderbank, D. M.~J., Gauduchon, P., and
  T\o{}nnesen-Friedman, C.~W.}
\newblock Hamiltonian 2-forms in {K}\"ahler geometry {III}, extremal metrics
  and stability.
\newblock {\em Invent. Math. 173}, 3 (2008), 547--601.

\bibitem{AP06}
{\sc Arezzo, C., and Pacard, F.}
\newblock Blowing up and desingularizing constant scalar curvature {K}\"ahler
  manifolds.
\newblock {\em Acta Math. 196}, 2 (2006), 179--228.

\bibitem{AT03}
{\sc Arezzo, C., and Tian, G.}
\newblock Infinite geodesic rays in the space of {K}\"ahler potentials.
\newblock {\em Ann. Sc. Norm. Super. Pisa Cl. Sci. (5) 2}, 4 (2003), 617--630.

\bibitem{AB83}
{\sc Atiyah, M.~F., and Bott, R.}
\newblock The {Y}ang-{M}ills equations over {R}iemann surfaces.
\newblock {\em Philos. Trans. Roy. Soc. London Ser. A 308\/} (1983), 523--615.

\bibitem{Aub78}
{\sc Aubin, T.}
\newblock {\'E}quations du type {M}onge-{A}mp\`ere sur les variet\'es
  k\"ahl\'eriennes compactes.
\newblock {\em Bull. Sci. Math. (2) 102}, 1 (1978), 63--95.

\bibitem{Aub84}
{\sc Aubin, T.}
\newblock R\'eduction de cas positif de l'\'equation de {M}onge-{A}mp\`ere sur
  les vari\'et\'es k\"ahl\'eriennes compactes \`a la d\'emonstration d'une
  in\'egalit\'e.
\newblock {\em J. Funct. Anal. 57}, 2 (1984), 143--153.

\bibitem{Bam16}
{\sc Bamler, R.}
\newblock {Convergence of Ricci flows with bounded scalar curvature}.
\newblock {\em arXiv:1603.05235\/}.

\bibitem{BKN89}
{\sc Bando, S., Kasue, A., and Nakajima, H.}
\newblock {On a construction of coordinates at infinity on manifolds with fast
  curvature decay and maximal volume growth}.
\newblock {\em Invent. Math. 97}, 2 (1989), 313--349.

\bibitem{BM85}
{\sc Bando, S., and Mabuchi, T.}
\newblock Uniqueness of {E}instein {K}\"ahler metrics modulo connected group
  actions.
\newblock In {\em Algebraic geometry, Sendai\/} (1985), vol.~10 of {\em Adv.
  Stud. Pure Math.}, pp.~11--40.

\bibitem{BT76}
{\sc Bedford, E., and Taylor, B.~A.}
\newblock The {D}irichlet problem for a complex {M}onge-{A}mp\`ere equation.
\newblock {\em Invent. Math. 37}, 1 (1976), 1--44.

\bibitem{Ber14}
{\sc Berman, R.}
\newblock {On the optimal regularity of weak geodesics in the space of metrics
  on a polarized manifold}.
\newblock {\em arXiv:1405.6482\/}.

\bibitem{Ber12}
{\sc Berman, R.}
\newblock K-polystability of {Q}-{F}ano varieties admitting
  {K}\"ahler-{E}instein metrics.
\newblock {\em Invent. Math. 203}, 3 (2016), 973--1025.

\bibitem{BB14}
{\sc Berman, R., and Berndtsson, B.}
\newblock Convexity of the {K}-energy on the space of {K}\"ahler metrics and
  uniqueness of extremal metrics.
\newblock {\em J. Amer. Math. Soc. 30}, 4 (2017), 1165--1196.

\bibitem{BBJ15}
{\sc Berman, R., Boucksom, S., and Jonsson, M.}
\newblock A variational approach to the {Y}au-{T}ian-{D}onaldson conjecture.
\newblock {\em arXiv:1509.04561\/}.

\bibitem{BDL16}
{\sc Berman, R., Darvas, T., and Lu, C.~H.}
\newblock {Regularity of weak minimizers of the K-energy and applications to
  properness and K-stability}.
\newblock {\em arXiv:1602.03114\/}.

\bibitem{BG14}
{\sc Berman, R., and Guenancia, H.}
\newblock {K\"ahler-Einstein metrics on stable varieties and log canonical
  pairs}.
\newblock {\em Geom. and Func. Anal. 24}, 6 (2014), 1683--1730.

\bibitem{Ber09}
{\sc Berndtsson, B.}
\newblock Curvature of vector bundles associated to holomorphic fibrations.
\newblock {\em Ann. of Math. (2) 169}, 2 (2009), 531--560.

\bibitem{Ber13}
{\sc Berndtsson, B.}
\newblock A {B}runn-{M}inkowski type inequality for {F}ano manifolds and some
  uniqueness theorems in {K}\"ahler geometry.
\newblock {\em Invent. Math. 200}, 1 (2015), 149--200.

\bibitem{Bl12}
{\sc B\l{}ocki, Z.}
\newblock {\em {On geodesics in the space of K\"ahler metrics}}, vol.~21 of
  {\em Adv. Lect. Math. (ALM)}.
\newblock Int. Press, Somerville, MA, 2012, pp.~3--19.

\bibitem{BFJ08}
{\sc Boucksom, S., Favre, C., and Jonsson, M.}
\newblock {Valuations and plurisubharmonic singularities}.
\newblock {\em Publ. RIMS 44\/} (2008), 449--494.

\bibitem{BHJ15}
{\sc Boucksom, S., Hisamoto, T., and Jonsson, M.}
\newblock {Uniform K-stability, Duistermaat-Heckman measures and singularities
  of pairs}.
\newblock {\em arXiv:1504.06568\/}.

\bibitem{BG08}
{\sc Boyer, C., and Galicki, K.}
\newblock {\em {Sasakian geometry}}.
\newblock {Oxford Mathematical Monographs}. OUP, 2008.

\bibitem{BrThesis}
{\sc Br\"onnle, T.~A.}
\newblock {\em Deformation constructions of extremal metrics}.
\newblock Phd Thesis, Imperial College London.

\bibitem{Cal54}
{\sc Calabi, E.}
\newblock The space of {K}\"ahler metrics.
\newblock {\em Proc. Internat. Congress Math. Amsterdm}, 2 (1954), 206--207.

\bibitem{Cal82}
{\sc Calabi, E.}
\newblock Extremal {K}\"ahler metrics.
\newblock In {\em Seminar on Differential Geometry}, S.~T. Yau, Ed. Princeton,
  1982.

\bibitem{CGP13}
{\sc Campana, F., Guenancia, H., and P{\u a}un, M.}
\newblock {Metrics with cone singularities along normal crossing divisors and
  holomorphic tensor fields}.
\newblock {\em Ann. Sci. {\'E}cole Norm. Sup. (4) 46\/} (2013), 879--916.

\bibitem{Cao85}
{\sc Cao, H.~D.}
\newblock Deformation of {K}\"ahler metrics to {K}\"ahler-{E}instein metrics on
  compact {K}\"ahler manifolds.
\newblock {\em Invent. Math. 81}, 2 (1985), 359--372.

\bibitem{Cat99}
{\sc Catlin, D.}
\newblock The {B}ergman kernel and a theorem of {T}ian.
\newblock In {\em Analysis and geometry in several complex variables (Katata,
  1997)\/} (1999), Trends Math., Birkh\"auser Boston Inc., pp.~1--23.

\bibitem{Ch01}
{\sc Cheeger, J.}
\newblock {\em Degeneration of {R}iemannian metrics under {R}icci curvature
  bounds}.
\newblock Lezioni Fermiane. Scuola Normale Superiore, Pisa, 2001.

\bibitem{CC97}
{\sc Cheeger, J., and Colding, T.}
\newblock On the structure of spaces with {R}icci curvature bounded below. {I}.
\newblock {\em J. Differential Geom. 46}, 3 (1997), 406--480.

\bibitem{CC00}
{\sc Cheeger, J., and Colding, T.}
\newblock On the structure of spaces with {R}icci curvature bounded below.
  {II}.
\newblock {\em J. Differential Geom. 54}, 1 (2000), 13--35.

\bibitem{CC00_1}
{\sc Cheeger, J., and Colding, T.}
\newblock On the structure of spaces with {R}icci curvature bounded below.
  {III}.
\newblock {\em J. Differential Geom. 54}, 1 (2000), 37--74.

\bibitem{CE08}
{\sc Cheeger, J., and Ebin, D.~G.}
\newblock {\em {Comparison theorems in Riemannian geometry}}.
\newblock AMS Chelsea Publishing, Providence, RI, 2008.

\bibitem{CN11}
{\sc Cheeger, J., and Naber, A.}
\newblock Lower bounds on {R}icci curvature and quantitative behavior of
  singular sets.
\newblock {\em Invent. Math. 191}, 2 (2013), 321--339.

\bibitem{CN14}
{\sc Cheeger, J., and Naber, A.}
\newblock Regularity of {E}instein manifolds and the codimension 4 conjecture.
\newblock {\em Ann. of Math. (2) 182}, 3 (2015), 1093--1165.

\bibitem{CP02}
{\sc Cheltsov, I., and Park, J.}
\newblock Log canonical thresholds and generalized {E}ckardt points.
\newblock {\em Sb. Math. 193}, 5--6 (2002), 779--789.

\bibitem{CDS12}
{\sc Chen, X., Donaldson, S., and Sun, S.}
\newblock K{\"a}hler-{E}instein metrics and stability.
\newblock {\em Int. Math. Res. Not. IMRN}, 8 (2014), 2119--2125.

\bibitem{CDS13_1}
{\sc Chen, X., Donaldson, S., and Sun, S.}
\newblock K{\"a}hler-{E}instein metrics on {F}ano manifolds. {I}:
  {A}pproximation of metrics with cone singularities.
\newblock {\em J. Amer. Math. Soc. 28}, 1 (2015), 183--197.

\bibitem{CDS13_2}
{\sc Chen, X., Donaldson, S., and Sun, S.}
\newblock K{\"a}hler-{E}instein metrics on {F}ano manifolds. {II}: {L}imits
  with cone angle less than {$2\pi$}.
\newblock {\em J. Amer. Math. Soc. 28}, 1 (2015), 199--234.

\bibitem{CDS13_3}
{\sc Chen, X., Donaldson, S., and Sun, S.}
\newblock K{\"a}hler-{E}instein metrics on {F}ano manifolds. {III}: {L}imits as
  cone angle approaches {$2\pi$} and completion of the main proof.
\newblock {\em J. Amer. Math. Soc. 28}, 1 (2015), 235--278.

\bibitem{CS12}
{\sc Chen, X., and Sun, S.}
\newblock {Space of K\"ahler metrics V -- K\"ahler quantization}.
\newblock In {\em Metric and differential geometry}, vol.~297 of {\em Progr.
  Math.} Birkh\"auser Verlag, Basel, 2012, pp.~19--41.

\bibitem{CSW15}
{\sc Chen, X., Sun, S., and Wang, B.}
\newblock {K\"ahler-Ricci flow, K\"ahler-Einstein metric, and K-stability}.
\newblock {\em arXiv:1508.04397\/}.

\bibitem{Chen00_1}
{\sc Chen, X.~X.}
\newblock The space of {K}\"ahler metrics.
\newblock {\em J. Differential Geom. 56}, 2 (2000), 189--234.

\bibitem{CT08}
{\sc Chen, X.~X., and Tang, Y.}
\newblock {Test configuration and geodesic rays}.
\newblock {\em Ast\'erisque}, 321 (2008), 139--167.

\bibitem{CW14}
{\sc Chen, X.~X., and Wang, B.}
\newblock {Space of Ricci flows (II)}.
\newblock {\em arXiv:1405.6797\/}.

\bibitem{CTW16}
{\sc Chu, J., Tosatti, V., and Weinkove, B.}
\newblock {On the $C^{1,1}$ regularity of geodesics in the space of K\"ahler
  metrics}.
\newblock {\em arXiv:1611.02390\/}.

\bibitem{Col97}
{\sc Colding, T.~H.}
\newblock Ricci curvature and volume convergence.
\newblock {\em Ann. of Math. (2) 145}, 3 (1997), 477--501.

\bibitem{CSz2}
{\sc Collins, T., and Sz\'ekelyhidi, G.}
\newblock {S}asaki-{E}instein metrics and {K}-stability.
\newblock {\em arXiv:1512.07213\/}.

\bibitem{Dar12}
{\sc Darvas, T.}
\newblock Morse theory and geodesics in the space of {K}\"ahler metrics.
\newblock {\em Proc. Amer. Math. Soc. 142}, 8 (2014), 2775--2782.

\bibitem{Dar15}
{\sc Darvas, T.}
\newblock {The Mabuchi geometry of finite energy classes}.
\newblock {\em Adv. Math. 285\/} (2015), 182--219.

\bibitem{DR15}
{\sc Darvas, T., and Rubinstein, Y.~A.}
\newblock {Tian's properness conjectures and Finsler geometry of the space of
  K\"ahler metrics}.
\newblock {\em J. Amer. Math. Soc. 30}, 2 (2017), 347--387.

\bibitem{DSz15}
{\sc Datar, V., and Sz\'ekelyhidi, G.}
\newblock {K}\"ahler-{E}instein metrics along the smooth continuity method.
\newblock {\em Geom. and Func. Anal. 26}, 4 (2016), 975--1010.

\bibitem{Del16}
{\sc Delcroix, T.}
\newblock {K}-stability of {F}ano spherical varieties.
\newblock {\em arXiv:1608.01852\/}.

\bibitem{Del15}
{\sc Delcroix, T.}
\newblock {K\"ahler-Einstein metrics on group compactifications}.
\newblock {\em Geom. and Func. Anal. 27}, 1 (2017), 78--129.

\bibitem{DK01}
{\sc Demailly, J.-P., and Koll\'ar, J.}
\newblock Semi-continuity of complex singularity exponents and
  {K}\"ahler-{E}instein metrics on {F}ano orbifolds.
\newblock {\em Ann. Sci. \'Ec. Norm. Sup\'er. (4) 34\/} (2001), 525--556.

\bibitem{Der16}
{\sc Dervan, R.}
\newblock {Relative K-stability for K\"ahler manifolds}.
\newblock {\em arXiv:1611.00569\/}.

\bibitem{Der15}
{\sc Dervan, R.}
\newblock {Alpha invariants and K-stability for general polarizations of Fano
  varieties}.
\newblock {\em Int. Math. Res. Not. IMRN}, 16 (2015), 7162--7189.

\bibitem{Der14}
{\sc Dervan, R.}
\newblock Uniform stability of twisted constant scalar curvature {K}\"ahler
  metrics.
\newblock {\em Int. Math. Res. Not. 15\/} (2016), 4728--4783.

\bibitem{DR16}
{\sc Dervan, R., and Ross, J.}
\newblock {K-stability for K\"ahler manifolds}.
\newblock {\em arXiv:1602.08983\/}.

\bibitem{Ding88}
{\sc Ding, W.~Y.}
\newblock Remarks on the existence problem of positive {K}\"ahler-{E}instein
  metrics.
\newblock {\em Math. Ann. 282}, 3 (1988), 463--471.

\bibitem{DT92}
{\sc Ding, W.~Y., and Tian, G.}
\newblock K\"ahler-{E}instein metrics and the generalized {F}utaki invariant.
\newblock {\em Invent. Math. 110}, 2 (1992), 315--335.

\bibitem{DS15}
{\sc Donaldson, S., and Sun, S.}
\newblock Gromov-{H}ausdorff limits of {K}\"ahler manifolds and algebraic
  geometry, {II}.
\newblock {\em arXiv:1507.05082\/}.

\bibitem{DS12}
{\sc Donaldson, S., and Sun, S.}
\newblock Gromov-{H}ausdorff limits of {K}{\"a}hler manifolds and algebraic
  geometry.
\newblock {\em Acta Math. 213}, 1 (2014), 63--106.

\bibitem{Don15}
{\sc Donaldson, S.~K.}
\newblock {The Ding functional, Berndtsson convexity and moment maps}.
\newblock {\em arXiv:1503.05173\/}.

\bibitem{Don83_1}
{\sc Donaldson, S.~K.}
\newblock {A new proof of a theorem of Narasimhan and Seshadri}.
\newblock {\em J. Differential Geom. 18}, 2 (1983), 269--277.

\bibitem{Don85}
{\sc Donaldson, S.~K.}
\newblock Anti self-dual yang-mills connections over complex algebraic surfaces
  and stable vector bundles.
\newblock {\em Proc. London Math. Soc. 50\/} (1985), 1--26.

\bibitem{Don87}
{\sc Donaldson, S.~K.}
\newblock {Infinite determinants, stable bundles and curvature}.
\newblock {\em Duke Math. J. 54}, 1 (1987), 231--247.

\bibitem{Don97}
{\sc Donaldson, S.~K.}
\newblock Remarks on gauge theory, complex geometry and four-manifold topology.
\newblock In {\em Fields Medallists' Lectures}, Atiyah and Iagolnitzer, Eds.
  World Scientific, 1997, pp.~384--403.

\bibitem{Don99_1}
{\sc Donaldson, S.~K.}
\newblock Symmetric spaces, {K}\"ahler geometry and {H}amiltonian dynamics.
\newblock In {\em Northern California Symplectic Geometry Seminar}, vol.~196 of
  {\em Amer. Math. Soc. Transl. Ser. 2}. Amer. Math. Soc., Providence, RI,
  1999, pp.~13--33.

\bibitem{Don01_1}
{\sc Donaldson, S.~K.}
\newblock {Planck's constant in complex and almost-complex geometry}.
\newblock In {\em {XIIIth International Congress on Mathematical Physics
  (London, 2000)}}. Int. Press, Boston, MA, 2001, pp.~63--72.

\bibitem{Don01}
{\sc Donaldson, S.~K.}
\newblock Scalar curvature and projective embeddings, {I}.
\newblock {\em J. Differential Geom. 59\/} (2001), 479--522.

\bibitem{Don02_1}
{\sc Donaldson, S.~K.}
\newblock {Holomorphic discs and the complex Monge-Amp\`ere equation}.
\newblock {\em J. Symplectic Geom. 1}, 2 (2002), 171--196.

\bibitem{Don02}
{\sc Donaldson, S.~K.}
\newblock Scalar curvature and stability of toric varieties.
\newblock {\em J. Differential Geom. 62\/} (2002), 289--349.

\bibitem{Don05}
{\sc Donaldson, S.~K.}
\newblock Lower bounds on the {C}alabi functional.
\newblock {\em J. Differential Geom. 70}, 3 (2005), 453--472.

\bibitem{Don04}
{\sc Donaldson, S.~K.}
\newblock Scalar curvature and projective embeddings, {II}.
\newblock {\em Q. J. Math. 56}, 3 (2005), 345--356.

\bibitem{Don09}
{\sc Donaldson, S.~K.}
\newblock Discussion of the {K}\"ahler-{E}instein problem.
\newblock {\em {\tt
  http://www2.imperial.ac.uk/{\textasciitilde}skdona/KENOTES.PDF}\/} (2009).

\bibitem{Don09_1}
{\sc Donaldson, S.~K.}
\newblock {Some numerical results in complex differential geometry}.
\newblock {\em Pure Appl. Math. Q. 5}, 2 (2009), 571--618.

\bibitem{Don12}
{\sc Donaldson, S.~K.}
\newblock {K\"ahler metrics with cone singularities along a divisor}.
\newblock In {\em Essays in mathematics and its applications}. Springer,
  Heilderlberg, 2012, pp.~49--79.

\bibitem{Don10}
{\sc Donaldson, S.~K.}
\newblock {\em Stability, birational transformations and the
  {K}\"ahler-{E}instein problem}, vol.~17 of {\em Surv. Differ. Geom.}
\newblock Int. Press, Boston, MA, 2012, pp.~203--228.

\bibitem{Don11}
{\sc Donaldson, S.~K.}
\newblock b-stability and blow-ups.
\newblock {\em Proc. Edinb. Math. Soc. (2) 57}, 1 (2014), 125--137.

\bibitem{DK90}
{\sc Donaldson, S.~K., and Kronheimer, P.~B.}
\newblock {\em The Geometry of Four-Manifolds}.
\newblock OUP, 1990.

\bibitem{Eys16}
{\sc Eyssidieux, P.}
\newblock {M\`etriques de K\"ahler-Einstein sur les vari\'et\'es de Fano
  [d'apr\`es Chen-Donaldson-Sun et Tian]}.
\newblock {\em Ast\'erisque}, {380, S\'eminaire Bourbaki. Vol. 2014/2015}
  (2016).

\bibitem{EGZ09}
{\sc Eyssidieux, P., Guedj, V., and Zeriahi, A.}
\newblock Singular {K}\"ahler-{E}instein metrics.
\newblock {\em J. Amer. Math. Soc. 22}, 3 (2009), 607--639.

\bibitem{Fine10}
{\sc Fine, J.}
\newblock {Calabi flow and projective embeddings. With an appendix by Kefeng
  Liu and Xiaonan Ma}.
\newblock {\em J. Differential Geom. 84}, 3 (2010), 489--523.

\bibitem{Fine12}
{\sc Fine, J.}
\newblock {Quantisation and the Hessian of Mabuchi energy}.
\newblock {\em Duke Math. J. 161}, 14 (2012), 2753--2798.

\bibitem{Fuj92}
{\sc Fujiki, A.}
\newblock Moduli space of polarized algebraic manifolds and {K}\"ahler metrics
  [translation of {S}\^ugaku {\bf 42} (1990), no. 3, 231--243;].
\newblock {\em Sugaku Expositions 5}, 2 (1992), 173--191.

\bibitem{Fu16}
{\sc Fujita, K.}
\newblock {K-stability of Fano manifolds with not small alpha invariants}.
\newblock {\em arXiv:1606.08261\/}.

\bibitem{Fuj15_1}
{\sc Fujita, K.}
\newblock {Optimal bounds for the volumes of K\"ahler-Einstein Fano manifolds}.
\newblock {\em arXiv:1508.04578\/}.

\bibitem{Fut83}
{\sc Futaki, A.}
\newblock An obstruction to the existence of {E}instein-{K}\"ahler metrics.
\newblock {\em Invent. Math. 73\/} (1983), 437--443.

\bibitem{Gau15}
{\sc Gauduchon, P.}
\newblock {\em {Calabi's extremal K\"ahler metrics: An elementary
  introduction}}.

\bibitem{GMSY}
{\sc Gauntlett, J.~P., Martelli, D., Sparks, J., and Yau, S.-T.}
\newblock Obstructions to the existence of {S}asaki-{E}instein metrics.
\newblock {\em Comm. Math. Phys. 273}, 3 (2007), 803--827.

\bibitem{Gro07}
{\sc Gromov, M.}
\newblock {\em Metric structures for {R}iemannian and non-{R}iemannian spaces}.
\newblock Birkh\"auser Boston Inc., 2007.

\bibitem{GT16}
{\sc Guenancia, H., and Taji, B.}
\newblock {Orbifold stability and Miyaoka-Yau inequality for minimal pairs}.
\newblock {\em arXiv:1611.0598\/}.

\bibitem{Ham82}
{\sc Hamilton, R.}
\newblock Three-manifolds with positive {R}icci curvature.
\newblock {\em J. Differential Geom. 17}, 2 (1982), 255--306.

\bibitem{HNpreprint}
{\sc Hein, H.-J., and Naber, A.}
\newblock {Isolated Einstein singularities with singular tangent cones}.
\newblock {\em in preparation\/}.

\bibitem{HS16}
{\sc Hein, H.-J., and Sun, S.}
\newblock {Calabi-Yau manifolds with isolated conical singularities}.
\newblock {\em arXiv:1607.02940\/}.

\bibitem{IS15}
{\sc Ilten, N., and S\"uss, H.}
\newblock K-stability for varieties with torus action of complexity one.
\newblock {\em Duke Math. J. 166}, 1 (2017), 177--204.

\bibitem{JMR16}
{\sc Jeffres, T., Mazzeo, R., and Rubinstein, Y.~A.}
\newblock {K\"ahler-Einstein metrics with edge singularities}.
\newblock {\em Ann. of Math. (2) 183}, 1 (2016), 95--176.

\bibitem{Jiang13}
{\sc Jiang, W.}
\newblock Bergman kernel along the {K}\"ahler {R}icci flow and {T}ian's
  conjecture.
\newblock {\em J. Reine Angew. Math. 717\/} (2016), 195--226.

\bibitem{Kah33}
{\sc K\"ahler, E.}
\newblock {\"Uber eine bemerkenswerte Hermitesche Metrik}.
\newblock {\em Abh. Math. Sem. Univ. Hamburg 9}, 1 (1933), 173--186.

\bibitem{KW74}
{\sc Kazdan, J.~L., and Warner, F.~W.}
\newblock Curvature functions for compact 2-manifolds.
\newblock {\em Ann. of Math. 99\/} (1974), 14--47.

\bibitem{KN79}
{\sc Kempf, G., and Ness, L.}
\newblock The length of vectors in representation spaces.
\newblock In {\em Algebraic geometry (Proc. Summer Meeting, Univ. Copenhagen,
  Copenhagen, 1978)}, vol.~732 of {\em Lecture Notes in Math.} Springer,
  Berlin, 1979, pp.~233--243.

\bibitem{Kol13}
{\sc Koll\'ar, J.}
\newblock Moduli of varieties of general type.
\newblock In {\em Handbook of moduli. {Vol. II}}. 2013, pp.~131--157.

\bibitem{KSB88}
{\sc Koll\'ar, J., and Shepherd-Barron, N.~I.}
\newblock Threefolds and deformations of surface singularities.
\newblock {\em Invent. Math. 91}, 2 (1988), 299--338.

\bibitem{Kol98}
{\sc Ko\l{}odziej, S.}
\newblock The complex {M}onge-{A}mp\`ere equation.
\newblock {\em Acta Math. 180}, 1 (1998), 69--117.

\bibitem{LV11}
{\sc Lempert, L., and Vivas, L.}
\newblock Geodesics in the space of {K}\"ahler metrics.
\newblock {\em Duke Math. J. 162}, 7 (2013), 1369--1381.

\bibitem{Li15}
{\sc Li, C.}
\newblock {K-semistability is equivariant volume minimization}.
\newblock {\em arXiv:1512.07205\/}.

\bibitem{Li11_2}
{\sc Li, C.}
\newblock {Constant scalar curvature K\"ahler metric obtains the minimum of
  K-energy}.
\newblock {\em Internat. Math. Res. Notices}, 9 (2011), 2161--2175.

\bibitem{LL16}
{\sc Li, C., and Liu, Y.}
\newblock {K\"ahler-Einstein metrics and volume minimization}.
\newblock {\em arXiv:1602.05094\/}.

\bibitem{LWX16}
{\sc Li, C., Wang, X., and Xu, C.}
\newblock {On proper moduli space of smoothable K\"ahler-Einstein Fano
  varieties}.
\newblock {\em arXiv:1411.0761\/}.

\bibitem{LWX15}
{\sc Li, C., Wang, X., and Xu, C.}
\newblock {Quasi-projectivity of the moduli space of smooth K\"ahler-Einstein
  Fano manifolds}.
\newblock {\em arXiv:1502.06532\/}.

\bibitem{LX16}
{\sc Li, C., and Xu, C.}
\newblock {Stability of Valuations and Koll\'ar Components}.
\newblock {\em arXiv:1604.05398\/}.

\bibitem{LX17_1}
{\sc Li, C., and Xu, C.}
\newblock {Stability of valuations: higher rational rank}.
\newblock {\em arXiv:1707.05561\/}.

\bibitem{LX11}
{\sc Li, C., and Xu, C.}
\newblock Special test configurations and {K}-stability of {F}ano varieties.
\newblock {\em Ann. of Math. (2) 180}, 2 (2014), 197--232.

\bibitem{LX17}
{\sc Liu, Y., and Xu, C.}
\newblock {K-stability of cubic threefolds}.
\newblock {\em arXiv:1706.01933\/}.

\bibitem{Lu98}
{\sc Lu, Z.}
\newblock On the lower order terms of the asymptotic expansion of
  {T}ian-{Y}au-{Z}elditch.
\newblock {\em Amer. J. Math. 122}, 2 (1998), 235--273.

\bibitem{Lun73}
{\sc Luna, D.}
\newblock Slices \'etales.
\newblock In {\em Sur les groupes alg\'ebriques}, {Bull. Soc. Math. France,
  Paris, M\'emoire 33}. Soc. Math. France, Paris, 1973, pp.~81--105.

\bibitem{Luo98}
{\sc Luo, H.}
\newblock {Geometric criterion for Gieseker-Mumford stability of polarized
  manifolds}.
\newblock {\em J. Differential Geom. 49\/} (1998), 577--599.

\bibitem{Mab86}
{\sc Mabuchi, T.}
\newblock K-energy maps integrating {F}utaki invariants.
\newblock {\em Tohoku Math. J. 38}, 4 (1986), 575--593.

\bibitem{Mab87}
{\sc Mabuchi, T.}
\newblock Some symplectic geometry on compact {K}\"ahler manifolds.
\newblock {\em Osaka J. Math. 24\/} (1987), 227--252.

\bibitem{Mab04}
{\sc Mabuchi, T.}
\newblock An energy-theoretic approach to the{H}itchin-{K}obayashi
  correspondence for manifolds, {I}.
\newblock {\em Invent. Math. 159}, 2 (2005), 225--243.

\bibitem{Mat57}
{\sc Matsushima, Y.}
\newblock Sur la structure du groupe d'hom\'eomorphismes analytiques d'une
  certaine vari\'et\'e k\"ahl\'erienne.
\newblock {\em Nagoya Math. J. 11\/} (1957), 145--150.

\bibitem{MFK94}
{\sc Mumford, D., Fogarty, J., and Kirwan, F.}
\newblock {\em Geometric invariant theory}, third~ed., vol.~34 of {\em
  Ergebnisse der Mathematik und ihrer Grenzgebiete (2) [Results in Mathematics
  and Related Areas (2)]}.
\newblock Springer-Verlag, Berlin, 1994.

\bibitem{Nad90}
{\sc Nadel, A.~M.}
\newblock Multiplier ideal sheaves and {K}\"ahler-{E}instein metrics of
  positive scalar curvature.
\newblock {\em Ann. of Math. (2) 132}, 3 (1990), 549--596.

\bibitem{NS65}
{\sc Narasimhan, M.~S., and Seshadri, C.~S.}
\newblock Stable and unitary vector bundles on a compact {R}iemann surface.
\newblock {\em Ann. of Math. (2) 82\/} (1965), 540--567.

\bibitem{Oda13_1}
{\sc Odaka, Y.}
\newblock {A generalization of the Ross-Thomas slope theory}.
\newblock {\em Osaka J. Math. 50}, 1 (2013), 171--185.

\bibitem{Oda13}
{\sc Odaka, Y.}
\newblock {The GIT stability of polarized varieties via discrepancy}.
\newblock {\em Ann. of Math. (2) 177}, 2 (2013), 645--661.

\bibitem{Oda15}
{\sc Odaka, Y.}
\newblock {Compact moduli spaces of K\"ahler-Einstein Fano varieties}.
\newblock {\em Publ. Res. Inst. Math. Sci. 51}, 3 (2015), 549--565.

\bibitem{OS12}
{\sc Odaka, Y., and Sano, Y.}
\newblock {Alpha invariant and K-stability of Q-Fano varieties}.
\newblock {\em Adv. Math. 229}, 5 (2012), 2818--2834.

\bibitem{OSS16}
{\sc Odaka, Y., Spotti, C., and Sun, S.}
\newblock {Compact moduli spaces of Del Pezzo surfaces and K\"ahler-Einstein
  metrics}.
\newblock {\em J. Differential Geom. 102}, 1 (2016), 127--172.

\bibitem{Paul12}
{\sc Paul, S.~T.}
\newblock Hyperdiscriminant polytopes, {C}how polytopes, and {M}abuchi energy
  asymptotics.
\newblock {\em Ann. of Math. (2) 175}, 1 (2012), 255--296.

\bibitem{PT06_1}
{\sc Paul, S.~T., and Tian, G.}
\newblock {CM stability and the generalized Futaki invariant, I.}
\newblock {\em arXiv:0605278\/}.

\bibitem{Per02}
{\sc Perelman, G.}
\newblock The entropy formula for the {R}icci flow and its geometric
  applications.
\newblock {\em math.DG/0211159\/}.

\bibitem{PSS12}
{\sc Phong, D.~H., Song, J., and Sturm, J.}
\newblock Degenerations of {K}\"ahler-{R}icci solitons on {F}ano manifolds.
\newblock {\em Univ. Iagel. Acta Math. 52\/} (2015), 29--43.

\bibitem{PS08}
{\sc Phong, D.~H., and Sturm, J.}
\newblock Lectures on stability and constant scalar curvature.
\newblock In {\em Current Developments in Mathematics 2007}. International
  Press.

\bibitem{PS06}
{\sc Phong, D.~H., and Sturm, J.}
\newblock Test-configurations for {K}-stability and geodesic rays.
\newblock {\em Jour. Symplectic Geom. 5}, 2 (2007), 221--247.

\bibitem{PS10}
{\sc Phong, D.~H., and Sturm, J.}
\newblock Regularity of geodesic rays and {M}onge-{A}mp\`ere equations.
\newblock {\em Proc. Amer. Math. Soc. 138}, 10 (2010), 3637--3650.

\bibitem{RT06}
{\sc Ross, J., and Thomas, R.~P.}
\newblock An obstruction to the existence of constant scalar curvature
  {K}\"ahler metrics.
\newblock {\em J. Differential Geom. 72\/} (2006), 429--466.

\bibitem{RT04}
{\sc Ross, J., and Thomas, R.~P.}
\newblock A study of the {H}ilbert-{M}umford criterion for the stability of
  projective varieties.
\newblock {\em J. Algebraic Geom. 16}, 2 (2007), 201--255.

\bibitem{RN15}
{\sc Ross, J., and Witt~Nystr\"om, D.}
\newblock {Harmonic discs of solutions to the complex homogeneous
  Monge-Amp\`ere equation}.
\newblock {\em Publ. Math. Inst. Hautes \'Etudes Sci. 122\/} (2015), 315--335.

\bibitem{RZ10}
{\sc Rubinstein, Y.~A., and Zelditch, S.}
\newblock {Bergman approximations of harmonic maps into the space of K\"ahler
  metrics on toric varieties}.
\newblock {\em J. Symplectic Geom. 8}, 3 (2010), 239--265.

\bibitem{ST15}
{\sc Sano, Y., and Tipler, C.}
\newblock Extremal metrics and lower bound of the modified {K}-energy.
\newblock {\em J. Eur. Math. Soc. 17}, 9 (2015), 2289--2310.

\bibitem{Sem92}
{\sc Semmes, S.}
\newblock Complex {M}onge-{A}mp\`ere equations and symplectic manifolds.
\newblock {\em Amer. J. Math. 114\/} (1992), 495--550.

\bibitem{Sey16}
{\sc Seyyedali, R.}
\newblock {Relative Chow stability and extremal metrics}.
\newblock {\em arXiv:1610.07555\/}.

\bibitem{SZ99}
{\sc Shiffman, B., and Zelditch, S.}
\newblock Distribution of zeros of random and quantum chaotic sections of
  positive line bundles.
\newblock {\em Comm. Math. Phys. 200}, 3 (1999), 661--683.

\bibitem{SD16}
{\sc Sj\"ostr\"om~Dyrefelt, Z.}
\newblock {K-semistability of cscK manifolds with transcendental cohomology
  class}.
\newblock {\em arXiv:1601.07659\/}.

\bibitem{SZ10}
{\sc Song, J., and Zelditch, S.}
\newblock Bergman metrics and geodesics in the space of {K}\"ahler metrics on
  toric varieties.
\newblock {\em Anal. PDE 3}, 3 (2010), 295--358.

\bibitem{Sp09}
{\sc Sparks, J.}
\newblock New results in {S}asaki-{E}instein geometry.
\newblock In {\em Riemannian topology and geometric structures on manifolds},
  vol.~271 of {\em Progr. Math.} Birkh\"auser Boston Inc., 2009, pp.~161--184.

\bibitem{SS17}
{\sc Spotti, C., and Sun, S.}
\newblock {Explicit Gromov-Hausdorff compactifications of moduli spaces of
  K\"ahler-Einstein Fano manifolds}.
\newblock {\em arXiv:1705.00377\/}.

\bibitem{SSY14}
{\sc Spotti, C., Sun, S., and Yao, C.}
\newblock Existence and deformations of {K}ahler-{E}instein metrics on
  smoothable {Q}-{F}ano varieties.
\newblock {\em Duke Math. J. 165}, 16 (2016), 3043--3083.

\bibitem{Sto08}
{\sc Stoppa, J.}
\newblock K-stability of constant scalar curvature {K}\"ahler manifolds.
\newblock {\em Adv. Math. 221}, 4 (2009), 1397--1408.

\bibitem{SSz09}
{\sc Stoppa, J., and Sz\'ekelyhidi, G.}
\newblock Relative {K}-stability of extremal metrics.
\newblock {\em J. Eur. Math. Soc. 13}, 4 (2011), 899--909.

\bibitem{Sz17}
{\sc Sz\'ekelyhidi, G.}
\newblock {Degenerations of $\mathbf{C}^n$ and Calabi-Yau metrics}.
\newblock {\em arXiv:1706.00357\/}.

\bibitem{GSz04}
{\sc Sz\'ekelyhidi, G.}
\newblock Extremal metrics and ${K}$-stability.
\newblock {\em Bull. Lond. Math. Soc. 39}, 1 (2007), 76--84.

\bibitem{GSz08}
{\sc Sz\'ekelyhidi, G.}
\newblock The {K}\"ahler-{R}icci flow and {K}-polystability.
\newblock {\em Amer. J. Math. 132\/} (2010), 1077--1090.

\bibitem{GSz13}
{\sc Sz\'ekelyhidi, G.}
\newblock Filtrations and test-configurations, with an appendix by {S.
  Boucksom}.
\newblock {\em Math. Ann. 362}, 1-2 (2015), 451--484.

\bibitem{Sz13_1}
{\sc Sz\'ekelyhidi, G.}
\newblock The partial {$C^0$}-estimate along the continuity method.
\newblock {\em J. Amer. Math. Soc. 29}, 2 (2016), 537--560.

\bibitem{Thomas06}
{\sc Thomas, R.~P.}
\newblock Notes on {GIT} and symplectic reduction for bundles and varieties.
\newblock In {\em Surveys in differential geometry. {V}ol. {X}}. Int. Press,
  Somerville, MA, 2006, pp.~221--273.

\bibitem{Tian87}
{\sc Tian, G.}
\newblock On {K}\"ahler-{E}instein metrics on certain {K}\"ahler manifolds with
  $c_1({M})>0$.
\newblock {\em Invent. Math. 89\/} (1987), 225--246.

\bibitem{Tian90_1}
{\sc Tian, G.}
\newblock On a set of polarized {K}\"ahler metrics on algebraic manifolds.
\newblock {\em J. Differential Geom. 32}, 1 (1990), 99--130.

\bibitem{Tian90}
{\sc Tian, G.}
\newblock On {C}alabi's conjecture for complex surfaces with positive first
  {C}hern class.
\newblock {\em Invent. Math. 101}, 1 (1990), 101--172.

\bibitem{Tian97}
{\sc Tian, G.}
\newblock K\"ahler-{E}instein metrics with positive scalar curvature.
\newblock {\em Invent. Math. 137\/} (1997), 1--37.

\bibitem{Tian00}
{\sc Tian, G.}
\newblock {\em Canonical metrics in {K}\"ahler geometry}.
\newblock Lectures in {M}athematics {ETH} {Z}\"urich. Birkh\"auser Verlag,
  Basel, 2000.

\bibitem{TZ12}
{\sc Tian, G., and Zhang, Z.}
\newblock Degeneration of {K}\"ahler-{R}icci solitons.
\newblock {\em Int. Math. Res. Not. IMRN}, 5 (2012), 957--985.

\bibitem{TZ13}
{\sc Tian, G., and Zhang, Z.}
\newblock Regularity of {K}\"ahler-{R}icci flows on {F}ano manifolds.
\newblock {\em Acta Math. 216}, 1 (2016), 127--176.

\bibitem{Tim11}
{\sc Timashev, D.~A.}
\newblock {\em {Homogeneous spaces and equivariant embeddings}}, vol.~138 of
  {\em {Encyclopaedia of Mathematical Sciences}}.
\newblock Springer, Heilderlberg, 2011.
\newblock {Invariant Theory and Algebraic Transformation Groups, 8}.

\bibitem{Uhl82_2}
{\sc Uhlenbeck, K.}
\newblock {Connections with $L^p$ bounds on curvature}.
\newblock {\em Comm. Math. Phys. 83\/} (1982), 31--42.

\bibitem{Uhl82_1}
{\sc Uhlenbeck, K.}
\newblock {Removable singularities in Yang-Mills fields}.
\newblock {\em Comm. Math. Phys. 83\/} (1982), 11--29.

\bibitem{UY86}
{\sc Uhlenbeck, K., and Yau, S.-T.}
\newblock On the existence of {H}ermitian-{Y}ang-{M}ills connections in stable
  vector bundles.
\newblock {\em Comm. Pure Appl. Math. 39\/} (1986), 257--293.

\bibitem{Wang12}
{\sc Wang, X.}
\newblock {Height and GIT weight}.
\newblock {\em Math. Res. Lett. 19}, 4 (2012), 909--926.

\bibitem{WZ04}
{\sc Wang, X.-J., and Zhu, X.}
\newblock K\"ahler-{R}icci solitons on toric manifolds with positive first
  {Chern} class.
\newblock {\em Adv. Math. 188\/} (2004), 87--103.

\bibitem{Ny10}
{\sc Witt~Nystr\"om, D.}
\newblock Test configurations and {O}kounkov bodies.
\newblock {\em Compos. Math. 148}, 6 (2012), 1736--1756.

\bibitem{Yau78}
{\sc Yau, S.-T.}
\newblock On the {R}icci curvature of a compact {K}\"ahler manifold and the
  complex {M}onge-{A}mp\`ere equation {I}.
\newblock {\em Comm. Pure Appl. Math. 31\/} (1978), 339--411.

\bibitem{Yau93}
{\sc Yau, S.-T.}
\newblock Open problems in geometry.
\newblock {\em Proc. Symposia Pure Math. 54\/} (1993), 1--28.

\bibitem{Zel98}
{\sc Zelditch, S.}
\newblock Szeg{\H o} kernel and a theorem of {T}ian.
\newblock {\em Int. Math. Res. Notices 6\/} (1998), 317--331.

\bibitem{Zh96}
{\sc Zhang, S.}
\newblock {Heights and reductions of semi-stable varieties}.
\newblock {\em Compos. Math. 104}, 1 (1996), 77--105.

\bibitem{Zhang16}
{\sc Zhang, Y.}
\newblock {Completion of the moduli space for polarized Calabi-Yau manifolds}.
\newblock {\em J. Differential Geom. 103}, 3 (2016), 521--544.

\end{thebibliography}

\end{document}